\documentclass[nohyperref]{article}

\usepackage{microtype}
\usepackage{graphicx}
\usepackage{subfigure}
\usepackage{booktabs} 

\usepackage{hyperref}

 \usepackage[accepted]{icml2022}

\usepackage{amsmath}
\usepackage{amssymb}
\usepackage{mathtools}
\usepackage{amsthm}
\usepackage{bm}
\usepackage{paralist}
\usepackage{multirow}

\usepackage[capitalize,noabbrev]{cleveref}

\theoremstyle{plain}
\newtheorem{theorem}{Theorem}[section]

\newtheorem{lemma}[theorem]{Lemma}

\theoremstyle{definition}

\newtheorem{assumption}[theorem]{Assumption}
\theoremstyle{remark}
\newtheorem{remark}[theorem]{Remark}

\newcommand{\bx}{\bm{x}}

\newcommand{\by}{\bm{y}}

\newcommand{\bp}{\bm{p}}

\newcommand{\bg}{\bm{g}}

\newcommand{\bw}{\bm{w}}

\newcommand{\bu}{\bm{u}}

\newcommand{\R}{\mathbb{R}}

\newcommand{\sO}{\mathcal{O}}

\usepackage[textsize=tiny]{todonotes}

\icmltitlerunning{Improve Single-Point Zeroth-Order Optimization Using High-Pass and Low-Pass Filters}

\begin{document}

\twocolumn[
\icmltitle{Improve  Single-Point Zeroth-Order Optimization \\Using High-Pass and Low-Pass Filters}



\icmlsetsymbol{equal}{*}

\begin{icmlauthorlist}
\icmlauthor{Xin Chen}{sch}
\icmlauthor{Yujie Tang}{sch}
\icmlauthor{Na Li}{sch}
\end{icmlauthorlist}

\icmlaffiliation{sch}{John A. Paulson School of Engineering and Applied Sciences,
Harvard University, MA, US}

\icmlcorrespondingauthor{Xin Chen}{chenxin2336@gmail.com}

\icmlkeywords{Machine Learning, ICML}

\vskip 0.3in
]



\printAffiliationsAndNotice{}  

\begin{abstract}
Single-point zeroth-order optimization (SZO) is useful in solving online black-box optimization and control problems in  time-varying environments, as it  queries  the function value only  once at each time step. However,  the vanilla SZO method is known to suffer from a large estimation variance and slow convergence, which seriously limits its practical application. In this work, we borrow 
the idea of
 high-pass and low-pass filters from extremum seeking  control (continuous-time version of SZO)  and develop a novel SZO method called HLF-SZO by integrating these filters. It turns out that the  high-pass filter coincides with the residual feedback  method, and the low-pass filter can be interpreted as the momentum method. As a result, the proposed HLF-SZO achieves a much smaller variance and much faster convergence than the vanilla SZO method, and empirically outperforms the residual-feedback SZO method,
 which are verified via extensive numerical experiments.
\end{abstract}

\section{Introduction}

This paper considers solving the generic  unconstrained optimization problem:
\begin{align}\label{eq:problem}
    \min_{\bx\in \R^d} \, f( \bx),
\end{align}
where $\bx\in \R^d$ is the decision variable and $f:\R^d\to \R$ is the objective function.
A straightforward solution  scheme is to apply the gradient descent method \cite{ruder2016overview}, e.g., $\bx\leftarrow \bx - \eta \nabla f(\bx)$ with the step size $\eta$.
 However, in many practical applications, the first-order information of function $f$, i.e., the gradient $\nabla f(\bx)$, may be unavailable or too expensive to procure, and
  one 
 can only access the zeroth-order  information,
i.e., function evaluations. 
To this end, 
zeroth-order (or derivative-free) methods \cite{nesterov2017random} are developed to solve   black-box   optimization and control problems, 
which essentially estimate the gradients using perturbed 
 function evaluations.
 Zeroth-order optimization (ZO) has attracted a great deal of recent attention and has
 been  used for a broad spectrum of applications, such as  reinforcement learning \cite{malik2019derivative,li2019distributed}, adversarial training \cite{chen2017zoo}, physical system control \cite{chen2021safe,chen2020model}, 
online sensor management \cite{liu2020primer}, etc.

\begin{table*}[t]
\caption{The iteration complexity of ZO methods for solving Lipschitz and smooth objective functions \eqref{eq:problem}. }
\label{tab:complexility}
\vskip 0.15in
\begin{center}
\begin{small}
\begin{sc}
\begin{tabular}{cccc}
\toprule
\multirow{2}{*}{Method} & \multirow{2}{*}{Literature} & \multicolumn{2}{c}{Complexity} \vspace{1pt}  \\\cline{3-4} \\[-0.7em]
    &  & $\quad\quad\ \ $ Convex $\quad\quad\ \ $  & Non-convex  \\
\midrule
Vanilla SZO   & \citet{gasnikov2017stochastic} & $\sO(d^{2}/\epsilon^{3}$)&  \textendash \\
Residual-feedback SZO & \citet{zhang2021new}& $\sO(d^{3}/\epsilon^{\frac{3}{2}})$ & $\sO(d^{3}/\epsilon^{\frac{3}{2}})$\\
\textbf{HLF-SZO}    & This work& $\sO(d^{\frac{3}{2}}/\epsilon^{\frac{3}{2}})$& $\sO(d^{\frac{3}{2}}/\epsilon^{\frac{3}{2}})$ \\
Two-point ZO    & \citet{nesterov2017random}& $\sO(d/\epsilon)$& $\sO(d/\epsilon)$         \\
\bottomrule \\[-0.7em]
\multicolumn{4}{l}{\small\normalfont *In the convex setting, the accuracy is measured by  $f(\overline{\bx}_T)\!-\!f(\bx^*)\leq \epsilon$, where $\bx^* \!\in\! \arg\min_{\bx\in\R^d} \! f(\bx)$ and $\overline{\bx}_T\!=\!\frac{1}{T}\sum_{k=1}^T\bx_k$.} \\
\multicolumn{4}{l}{\small\normalfont $\ $ In the non-convex setting, the accuracy is measured by $\frac{1}{T}\sum_{k=1}^T\|\nabla f(\bx_k)\|^2\leq \epsilon$. } \\
\end{tabular}
\end{sc}
\end{small}
\end{center}
\vskip -0.1in
\end{table*}

According to the number of queried
function evaluations at each iteration,  ZO methods can be categorized into two types: single-point and multi-point \cite{liu2020primer}. 
As suggested by the name, single-point ZO (SZO) \cite{flaxman2005online} only needs to query the function value once at each iteration, making it particularly suitable for online optimization and control problems. In contrast, multi-point ZO, such as two-point ZO \cite{shamir2017optimal}, requires two or more function evaluations in the same instantaneous time; this may not be 
practical when the environment is non-stationary or changed by the implementation of a function evaluation. For example, two-point ZO can not be applied to non-stationary reinforcement learning problems, because it needs two different policy evaluations in the same environment, which is impossible since the environment changes after each policy evaluation
\cite{zhang2021new}. 
Nevertheless, it is known that the vanilla SZO method \cite{flaxman2005online} suffers from  a large estimation variance and 
slow convergence, which seriously jeopardizes its practical application. There have been multiple studies \cite{saha2011improved,dekel2015bandit,gasnikov2017stochastic,zhang2021new} on improving the convergence rate of SZO. In particular,  the residual-feedback SZO method proposed in the recent work \citet{zhang2021new} achieves the
state-of-the-art performance   to date.

In the  field of control,  there is a continuous-time version of SZO  known as \emph{extremum seeking (ES) control} \cite{ariyur2003real,tan2010extremum}.  ES control is  a classic   adaptive  control technique that uses only output feedback to
steer a dynamical system to a state where the objective function attains an extremum.  ES control can also be 
adopted as a zeroth-order algorithm to solve static-map optimization problems \cite{PovedaNaLi2019,durr2013saddle,ye2016distributed}. 
In addition, although not essential for the ES system operation, a \emph{high-pass filter} and a \emph{low-pass filter} are usually integrated into the control loop, because these filters can significantly improve the transient behavior and mitigate oscillations
(see Section \ref{sec:pre:ES} for detailed explanations).
Despite their close connection, ES control and SZO have been mostly studied separately in the control and optimization communities. 
Then, a natural question to ask is 

 \emph{``Can we borrow some  tools   from extremum seeking control, such as the high-pass and low-pass filters, to improve the performance of single-point zeroth-order optimization?''}


\textbf{Contributions.} Motivated by this question,  we develop a 
 novel SZO method called \textbf{HLF-SZO} (\textbf{H}igh/\textbf{L}ow-pass \textbf{F}ilter \textbf{SZO}) by integrating a  high-pass filter and a low-pass filter into the vanilla SZO method. The main contributions of this paper are explained below:
 
 \setdefaultleftmargin{12pt}{0pt}{}{}{}{}
\begin{enumerate}
\item [1)]  We find that the integration of a high-pass filter can be interpreted as 
 the residual feedback  scheme proposed in the recent work \citet{zhang2021new},
which can greatly reduce the variance of SZO and lead to faster convergence. And the integration of a low-pass filter can be interpreted as the momentum (heavy-ball) optimization method \cite{polyak1964some,qian1999momentum}, which can further accelerate the convergence.


\item [2)]   We prove that the proposed  HLF-SZO method  achieves the iteration complexity of $\sO(d^{\frac{3}{2}}/\epsilon^{\frac{3}{2}})$ for Lipschitz and smooth objective functions in both convex and nonconvex cases. This iteration complexity 
is better than the  complexity $\sO(d^{2}/\epsilon^{3})$ of the vanilla SZO method \cite{gasnikov2017stochastic}, and has better dependency on the problem dimension $d$ compared with the  complexity $\sO(d^3/\epsilon^{\frac{3}{2}})$ of the residual-feedback SZO method \cite{zhang2021new},
although it is inferior to  the complexity $\sO(d/\epsilon)$ of two-point methods. See Table \ref{tab:complexility} for a
  summary of the best known iteration complexity of ZO methods.

\item [3)] Extensive numerical experiments 
  show that 
 the  proposed HLF-SZO method  exhibits a much smaller variance and much faster convergence than the vanilla SZO; it   empirically outperforms the residual-feedback SZO and has comparable performance to the two-point ZO method.
\end{enumerate}

\vspace{-10pt}
Moreover,
this paper  explores a new direction to improve ZO schemes by leveraging the connection between   ZO and continuous-time ES control, and we  exemplify the possibility that  useful tools from ES control, such as high-pass and low-pass filters,  can indeed help boost SZO. There is much more to study in this direction. 
 For example,  this paper only adopts the simplest high-pass and low-pass filters, while
higher-order filters or compensators that are used in ES control to  enhance  stability \cite{krstic2000performance} may also be applied to further improve the performance of ZO.

\textbf{Other Related Work on ZO.} The work \cite{novitskii2021improved} leverages the higher-order smoothness and proves an improved complexity bound of SZO for solving stochastic convex optimization problems.
References \cite{jongeneel2021small, jongeneel2021imaginary} propose new smoothed gradient approximation schemes by using complex analysis tools from numerical differentiation, which boost the convergence of ZO and address the numerical cancellation issue. In \citet{wang2021improved,wang2021model}, a gradient estimation method that uses    only the system measurements is developed based on   the simultaneous perturbation stochastic
approximation   algorithm and the complex-step gradient
approximation.
\citet{berahas2022theoretical} presents a theoretical and empirical comparison of several gradient approximation methods for derivative-free optimization, including finite difference, linear interpolation, Gaussian smoothing, etc.

\noindent{\bf Notations.} Let $\mathbb{B}_d\coloneqq \{\bx\in\mathbb{R}^d:\|\bx\|_2\leq 1\}$ denote the closed unit ball of dimension $d$, and let $\mathbb{S}_{d-1}\coloneqq \{\bx\in\mathbb{R}^d:\|\bx\|_2= 1\}$ denote the unit sphere.

\section{Preliminaries on ZO and ES Control}\label{sec:preliminary}

\subsection{Zeroth-Order Optimization}

The SZO method  estimates the gradient $\nabla f(\bx)$ with only one  query of the function value at each iteration step. Specifically, upon defining
\begin{equation}\label{eq:single_grad_est}
\mathsf{G}_f^{(1)}(\bx;r,\bu)
\coloneqq \frac{d}{r}f(\bx + r\bu )\bu,
\end{equation}
for $\bx,\bu\in\mathbb{R}^d$ and $r>0$, one can show that if $\bu$ is uniformly sampled from the unit sphere $\mathbb{S}_{d-1}$, the term $\mathsf{G}_f^{(1)}(\bx;r,\bu)$ is an unbiased estimator for the
gradient of a smoothed version of function $f$ \cite{flaxman2005online}:\footnote{
Another choice of the distribution to sample the random direction $\bu$ is the Gaussian distribution $\mathcal{N}(\bm{0},\bm{I}/d)$. We use $\mathrm{Unif}(\mathbb{S}_{d-1})$ in this paper for the simplicity of exposition, while the results can be also extended to the Gaussian distribution case.
}
\begin{align*}
    \mathbb{E}_{\bu\sim \mathrm{Unif}(\mathbb{S}_{d-1})}\big[ \mathsf{G}_f^{(1)}(\bx;r,\bu)\big] = \nabla f_r(\bx),
\end{align*}
where $f_r(\bx) := \mathbb{E}_{\bar{\bu}\sim \mathrm{Unif}(\mathbb{B}_d)}\big[ f(\bx+r \bar{\bu})\big]$.
By plugging the single-point gradient estimator $\mathsf{G}_f^{(1)}(\bx;r,\bu)$ into the gradient descent  iterations for solving \eqref{eq:problem}, we obtain the vanilla SZO method \eqref{eq:single:basic}:
\begin{align}\label{eq:single:basic}
 (\textbf{Vanilla SZO}):   \bx_{k+1} = \bx_k - \eta\cdot \frac{d}{r} f(\bx_k\!+\! r \bu_k) \bu_k,
\end{align}
where $\{\bu_k:k=0,1,2,\cdots\}$ are i.i.d. random directions 
 sampled from the uniform distribution $\mathrm{Unif}(\mathbb{S}_{d-1})$.
 Here, $\eta\!>\!0$ is the step size and the parameter $r\!>\!0$ is called  smoothing radius. One can show the convergence of~\eqref{eq:single:basic} by appropriately choosing the parameters  $\eta$ and  $r$ under some regular conditions \cite{gasnikov2017stochastic}.

However, the single-point gradient estimator~\eqref{eq:single_grad_est} generally suffers from a high estimation variance, which leads to slow convergence of the vanilla SZO method~\eqref{eq:single:basic}. To overcome this issue and achieve faster convergence, two-point ZO methods \cite{nesterov2017random} are developed as:
\begin{align}
    & (\textbf{Two-point ZO}):\nonumber \\  & \bx_{k+1}\! = \bx_k 
    \!-\! \eta \frac{d}{2r} \big(f(\bx_k\!+\! r \bu_k) -f(\bx_k\!-\! r \bu_k)\big)\bu_k,\label{eq:two}\\
   & \text{(or)}\ \ \bx_{k+1} = \bx_k  - \eta \frac{d}{r} \big(f(\bx_k\!+\! r \bu_k) -f(\bx_k)\big)\bu_k, \label{eq:two2}
\end{align}
where two function evaluations are needed in each iteration. Both of the two-point gradient estimators $\mathsf{G}_f^{(2)}(\bx;r,\bu)
\coloneqq \frac{d}{2r}(f(\bx+r\bu)-f(\bx-r\bu))\bu$ and $\Tilde{\mathsf{G}}_f^{(2)}(\bx;r,\bu)
\coloneqq \frac{d}{r}(f(\bx+r\bu)-f(\bx))\bu$ share the same expectation as  $\mathsf{G}_f^{(1)}(\bx;r,\bu)$, but  they generally have smaller variances and thus lead to faster convergence of the two-point ZO methods \eqref{eq:two}  \eqref{eq:two2}. 

\subsection{Extremum Seeking (ES) Control}\label{sec:pre:ES}

\begin{figure}
    \centering
  \includegraphics[scale = 0.35]{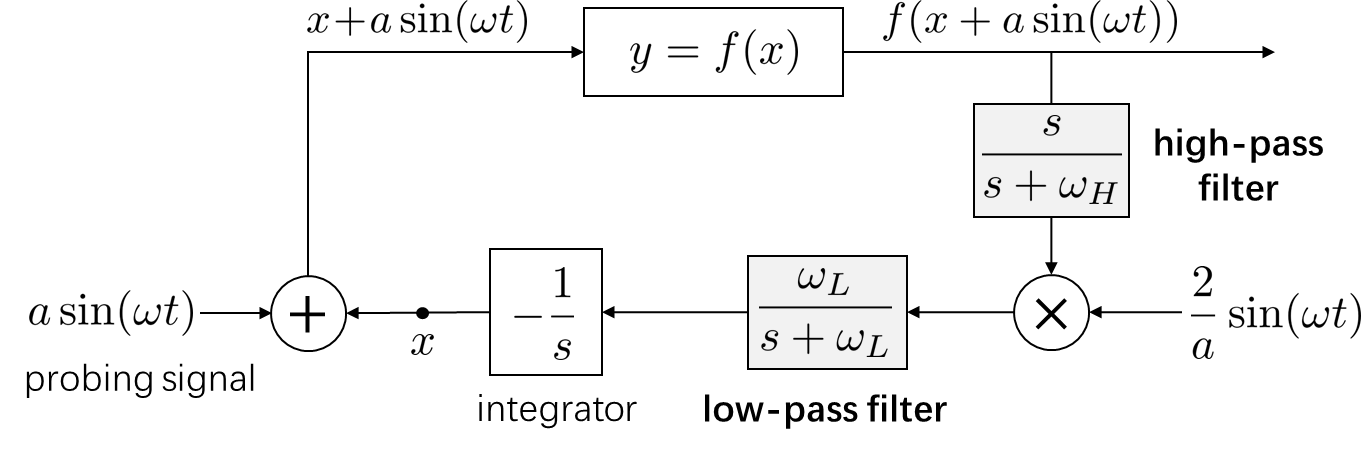}
    \caption{The block diagram of ES system with a high-pass filter $\frac{s}{s+\omega_H}$ and a low-pass filter 
    $\frac{\omega_L}{s+\omega_L}$ for solving $\min_x f(x)$.}
    \label{fig:ES}
\end{figure}

Extremum seeking control is a  type of model-free control that uses  only output feedback to steer a dynamical system to the extremum of an unknown function \cite{ariyur2003real}.  
A simple ES scheme for solving the scalar version of problem \eqref{eq:problem} is illustrated as Figure \ref{fig:ES}. Essentially, the ES system adopts a small sinusoidal  probing signal $a\sin(\omega t)$  as the perturbation to estimate the gradient $\nabla f(x)$, where positive parameters $a$ and $\omega$ are the amplitude and  frequency.

 To analyze the ES system shown in Figure~\ref{fig:ES}, we can first ignore 
the high-pass filter  and the low-pass filter, 
as they are not essential to the system operation. Following the block diagram in Figure~\ref{fig:ES}, we start from the state $x$ and add a sinusoidal probing signal $a\sin(\omega t)$  to it, and the resultant value $x+a\sin{\omega t}$ is taken as the input 
to the static map $y=f(x)$. Then, the  output $y=f(x+a\sin(\omega t))$  is  multiplied by $\frac{2}{a}\sin(\omega t)$ for correlation and the  loop is  closed after passing through the integrator $-\frac{1}{s}$. As a consequence, the 
ES system dynamics can be  formulated as 
\begin{align}\label{eq:es}
  & (\textbf{Extremum Seeking Dynamics}):\nonumber\\  &\qquad \dot{x} = - \frac{2}{a}f\big(x(t) + a\sin(\omega t) \big)\sin(\omega t).
\end{align}
The formulation of the ES dynamics \eqref{eq:es} is quite similar to the vanilla SZO method \eqref{eq:single:basic}, except that a sinusoidal probing signal is used instead of  random directions. 
 Besides,
the ES scheme can be easily extended to the multivariate case by using an appropriate probing signal vector with different frequencies
\cite{chen2021model,PovedaNaLi2019}.

 The rationale behind \eqref{eq:es} is that   with a small  $a$  and a large  $\omega$,  
 the ES dynamics \eqref{eq:es}  behaves approximately like the gradient descent flow  $\dot{x}=-\nabla f(x)$, which  steers  $x$ to a (local) minimum $x^*\in\arg\min_{x}f(x)$ under regular conditions \cite{chill2010gradient}. 
Specifically, when $\omega$ is large, 
 the ES dynamics \eqref{eq:es}  exhibits a timescale separation property, where the fast-time variation  results from the high-frequency sinusoidal signal $\sin(\omega t)$, while the slow-time variation induced by the integrator $-\frac{1}{s}$   dominates the evolution of $x$. As $a$ is small, 
we consider 
 the Taylor expansion  of the output $f(x+a\sin(\omega t))$ around $x$:
 \begin{align} \label{eq:taylor}
     f(x\!+\!a\sin(\omega t))\! =\! f(x) \!+\! a\sin(\omega t)\nabla f(x)\! +\! \sO(a^2).
 \end{align}
 Then to analyze the dominant slow-time variation,
 one can compute the average dynamics of \eqref{eq:es} to wash out the periodic fast-time variation, which is $\dot{x} = -\, h_{\mathrm{ave}}(x)$ with
\begin{align}
   & h_{\mathrm{ave}}(x) \coloneqq \, \frac{1}{T}\int_{0}^T \frac{2}{a}  f(x+a\sin(\omega t))\sin(\omega t) \, dt \nonumber\\
   & = \frac{1}{T}\!\int_{0}^T\! \frac{2}{a} f(x)\sin(\omega t) +   2\sin^2(\omega t) \nabla f(x) +   {\sO}(a)\, dt\nonumber \\
     & = \nabla f(x) + {\sO}(a), \label{eq:inte}
\end{align}
where $T=\frac{2\pi}{\omega}$. Equation~\eqref{eq:inte} reveals that  the average dynamics of \eqref{eq:es} is actually  the gradient descent flow plus a small perturbation term ${\sO}(a)$. See \citet{ariyur2003real,tan2010extremum} for more explanations.

\vspace{5pt}
\begin{remark}\label{remark:hlf}
(\textbf{Add High-Pass and Low-Pass Filters}).
As shown in Figure \ref{fig:ES}, 
 a  high-pass filter $\frac{s}{s+\omega_H}$ and a lower-pass filter $\frac{\omega_L}{s+\omega_L}$ are integrated into the ES system to mitigate oscillations and improve  the transient  behavior,
 where $\omega_H$ and $\omega_L$ are the cut-off frequency parameters.
  Intuitively,  the high-pass filter serves to wash out the constant and low-frequency components of the output signal $f(x+a\sin(\omega t))$. 
  As indicated in the integral \eqref{eq:inte}, the first  term $f(x)$ in the Taylor expansion \eqref{eq:taylor} 
  is useless and 
  causes large oscillations in the transient process, and it is 
  slowly-varying and thus can be removed by a high-pass filter. In terms of the low-pass filter, it works  to wash out high-frequency oscillations (induced by the sinusoidal probing signal) of the gradient estimation, because a clean and constant  gradient estimator is desirable before passing through the integrator. More explanations are provided in  \citet{tan2010extremum,tan2006non}.
\end{remark}

\section{Algorithm Design}

In this section, we first present the proposed SZO algorithms that incorporate  high-pass and low-pass filters, and then describe the derivation process.

\subsection{SZO Methods with High-Pass/Low-Pass Filters}

\subsubsection{Integrate a High-Pass Filter} 
We mimic how  ES control applies the high-pass  filter and integrate it into the vanilla SZO \eqref{eq:single:basic}, leading to the following new SZO method
\eqref{eq:single:high}: 
\begin{align} \label{eq:single:high}
  &\qquad (\textbf{HF-SZO}):\quad \qquad \forall k=1,2,\cdots\nonumber\\
  &\begin{dcases}
   z_k  \!= (1\!-\!\beta) z_{k-1} + f(\bx_k \!+\!r \bu_k) \!- f(\bx_{k-1} \!+\!r\bu_{k-1})\\
   \bx_{k+1} = \bx_k - \eta\cdot \frac{d}{r} z_k \bu_k, 
   \end{dcases}
\end{align}
where  $z_k\in \R$ is an intermediate variable 
and $\beta\geq 0$ is an adjustable  parameter. Particularly, when setting $\beta =1$, the HF-SZO 
 \eqref{eq:single:high} becomes the  residual-feedback SZO method proposed in   \citet{zhang2021new};  and \eqref{eq:single:high} reduces to the vanilla SZO \eqref{eq:single:basic} when   $\beta= 0$ and $z_0 = f(\bx_0+r\bu_0)$.

\subsubsection{Integrate a Low-Pass Filter}

Similarly, by  integrating a low-pass filter into the vanilla SZO \eqref{eq:single:basic}, we obtain a new SZO method \eqref{eq:single:low}:
\begin{align} \label{eq:single:low}
 &\quad (\textbf{LF-SZO}):\quad \qquad \forall k=1,2,\cdots\nonumber\\
&    \bx_{k+1}\! = \bx_{k} - {\eta}  \frac{d}{r} f(\bx_k\!+\! r \bu_k) \bu_k + \alpha (\bx_k\!-\! \bx_{k-1}),
\end{align}
where $\alpha\in[0,1]$ is an adjustable  parameter. 

Compared with \eqref{eq:single:basic}, the LF-SZO method  \eqref{eq:single:low} has an additional 
 ``momentum'' term $\alpha(\bx_k\!-\!\bx_{k-1})$ that is introduced by the low-pass filter.
When setting $\alpha=0$, the LF-SZO method \eqref{eq:single:low} reduces to the vanilla SZO method \eqref{eq:single:basic}.  

\subsubsection{Integrate Both a High-Pass Filter and a Low-Pass  Filter}

By integrating   both a high-pass filter and a low-pass filter into \eqref{eq:single:basic}, we develop the HLF-SZO method \eqref{eq:single:both}:
\begin{align} \label{eq:single:both}
&\quad (\textbf{HLF-SZO}):\quad \qquad \forall k=1,2,\cdots\nonumber\\
&\begin{dcases}
  z_k = (1\!-\!\beta) z_{k-1} + f(\bx_k \!+\!r \bu_k) \!-\! f(\bx_{k-1} \!+\! r\bu_{k-1})\\
   \bx_{k+1} = \bx_{k} - {\eta} \cdot \frac{d}{r} z_k \bu_k + \alpha (\bx_k- \bx_{k-1}),
 \end{dcases}
\end{align}
which is basically the combination of \eqref{eq:single:high} and \eqref{eq:single:low}.

Note that in the HLF-SZO \eqref{eq:single:both}, only one function evaluation $f(\bx_k +r \bu_k)$ is queried at each iteration $k$, while the value $f(\bx_{k-1} \!+\!r\bu_{k-1})$ is directly inherited from the last iteration, which is the same for the HF-SZO \eqref{eq:single:high}.
 As for the choices of adjustable parameters $\alpha$ and $\beta$, 
the theoretical analysis in Section \ref{sec:analysis} shows that $\beta =1$ is optimal for convergence, which is also validated by the 
numerical experiments conducted in Section \ref{sec:numerical}. Besides, a momentum parameter $\alpha=0.9$ or a similar value is  suggested  \cite{ruder2016overview} and there have been extensive studies \cite{polyak1964some,tao2021role,qian1999momentum} on the role and selection of the momentum parameter. See Section \ref{sec:analysis} for more discussions.

\vspace{5pt} 
\begin{remark}
High-pass and low-pass filters are classic tools in the field of control and signal processing. 
It is interesting to find their connections to optimization approaches. Specifically,
 the integration of a \emph{high-pass filter} coincides with the \emph{residual-feedback method} proposed in \citet{zhang2021new}, which can significantly reduce the estimation variance of SZO methods. This is also consistent with the function of a high-pass filter  explained in Remark \ref{remark:hlf}. In addition, the integration of a \emph{low-pass filter} can be interpreted as the \emph{momentum (heavy-ball) method} \cite{polyak1964some}, which can accelerate the convergence. These observations are further explained in Section \ref{sec:analysis} and are verified via extensive numerical experiments in Section \ref{sec:numerical}.
\end{remark}




\subsection{Detailed Derivation Process}

For the ES system illustrated in Figure \ref{fig:ES}, we denote $f$ as the output of the static map, and let
$z$ be the signal after passing $f$ through the high-pass filter $\frac{s}{s+\omega_H}$. Then we have 
\begin{align} \label{eq:highdyn}
  \mathcal{L}\{z\} = \frac{s}{s+\omega_H} \mathcal{L}\{f\}\ \Longleftrightarrow\ \dot{z} +\omega_H z = \dot{f},
\end{align}
where $\mathcal{L\{\cdot\}}$ denotes the Laplacian transform.   We discretize the continuous-time dynamics \eqref{eq:highdyn} under the time gap $\delta$:
\begin{align} \label{eq:high:med}
    &\frac{z_k -z_{k-1}}{\delta} + \omega_H z_{k-1} = \frac{f_k - f_{k-1}}{\delta}\nonumber \\
    &\qquad \Longleftrightarrow\ \  z_k = (1-\delta \omega_H) z_{k-1} + f_k - f_{k-1},
\end{align}
where $f_k\coloneqq f(\bx_k + r\bu_k)$. Essentially, $z_k$ can be regarded  as the refined value of $f_k$ after passing it through the high-pass filter. Denote $\beta: = \delta \omega_H\geq 0$. Then replacing $f_k$ by $z_k$ in the vanilla SZO \eqref{eq:single:basic} leads to the HF-SZO  method \eqref{eq:single:high}.


 As shown in Figure \ref{fig:ES}, a low-pass filter $ \frac{\omega_L}{s+\omega_L}$ is added  before passing the gradient estimation to the integrator. 
 We denote $\bg\in \R^d$ as the  gradient estimator and $\by\in\R^d$ as the refined signal after passing $\bg$ through the low-pass filter. Then we have the relation:
\begin{align}\label{eq:lowdyn}
    \mathcal{L}\{\by \} =\frac{\omega_L}{s+\omega_L} \mathcal{L}\{\bg \} \ \Longleftrightarrow\ \dot{\by} = \omega_L(-\by +\bg).
\end{align}
We discretize the resultant dynamics \eqref{eq:lowdyn} and the integrator dynamics $ \dot{\bx} = -\by$ under the time gap $\delta$, and obtain
\begin{subequations} \label{eq:low:med}
\begin{align}
   & \frac{\by_{k+1} -\by_k}{\delta} = \omega_L(-\by_k +\bg_k)\nonumber \\ &\qquad \Longleftrightarrow \ \ \by_{k+1} = (1-\delta \omega_L) \by_k  + \delta \omega_L \bg_k, \\
    & \frac{\bx_{k+1}-\bx_k}{\delta} = - \by_{k+1}   \Longleftrightarrow  \bx_{k+1} = \bx_k - \delta \by_{k+1}.
\end{align}
\end{subequations}
Denote ${\eta}\coloneqq \delta^2 \omega_L$ and $\alpha\coloneqq 1-\delta \omega_L$. After eliminating the variable $\by$ from \eqref{eq:low:med} and letting $\bg_k=  \frac{d}{r} f(\bx_k+ r \bu_k) \bu_k$ as
in the vanilla SZO \eqref{eq:single:basic}, 
 we obtain the LF-SZO \eqref{eq:single:low}.

When both a high-pass filter and a low-pass filter are applied,
we obtain the HLF-SZO method \eqref{eq:single:both}, which is basically
 the combination of \eqref{eq:high:med} and \eqref{eq:low:med}.

\section{Theoretic Analysis}\label{sec:analysis}

In this section, we analyze  the convergence of the proposed HLF-SZO method~\eqref{eq:single:both}. 
We make the following assumption throughout our analysis.
\begin{assumption}\label{ass:lip}
The function $f$ is $G$-Lipschitz and $L$-smooth, i.e., for all $\bx,\by\in\mathbb{R}^d$, we have
$$
|f(\bx)-f(\by)|\!\leq\! G\|\bx-\by\|,\,
\|\nabla\! f(\bx)-\nabla\! f(\by)\|\!\leq \! L\|\bx-\by\|.
$$
\end{assumption}

Note that Assumption \ref{ass:lip} is mainly for theoretical analysis, and the
proposed HLF-SZO works well for a wide range of problems that may not satisfy this assumption, as verified by the
numerical experiments in Section \ref{sec:numerical}.

We first consider the case when the objective function $f$ is convex, and the following theorem states the convergence properties of the HLF-SZO \eqref{eq:single:both}.

\begin{theorem}\label{theorem:main_convex}
(Convex Case). Let $\alpha\in[0,1)$, $\beta\in(0,2)$, and  the total number of iterations $T$ be given. Suppose that Assumption \ref{ass:lip} holds,
$f$ is convex and has a finite minimizer $\bx^\ast\in\mathbb{R}^d$.  Then by choosing
\begin{equation}\label{eq:condition_convergence}
\eta \!\leq\! \frac{(1\!-\!\alpha)(1\!-\!|1\!-\!\beta|)^2}{20Ld T^{\frac{1}{3}}},\, 
\frac{4\eta d G}{(1\!-\!|1\!-\!\beta|)(1\!-\!\alpha)}\!\leq\! r \!\leq\! \frac{G}{L T^{\frac{1}{3}}},
\end{equation}
the HLF-SZO method \eqref{eq:single:both}  achieves
\begin{align}
    &\mathbb{E}\left[f(\overline{\bx}_T)\right]-f(\bx^\ast)\nonumber\\
&\ \  \leq
\frac{3(1-\alpha)\|\bx_1-\bx^\ast\|^2}{4\eta T} +
\frac{3G^2}{2 L T^{2/3}}
+ \sO(\frac{d}{T}),
\end{align}
where $\overline{\bx}_T\coloneqq \frac{1}{T}\sum_{k=1}^T\bx_k$. Moreover, by letting $\eta$ achieve the equality in~\eqref{eq:condition_convergence}, we have
\begin{align}\label{eq:convergence_convex}
&\mathbb{E}\left[f(\overline{\bx}_T)\right]-f(\bx^\ast)\nonumber\\
&\ \ \leq
\frac{d}{T^{2/3}}
\left(
\frac{15 L \|\bx_1-\bx^\ast\|^2}{(1-|1-\beta|)^2}
+\frac{3G^2}{2Ld}\right)
+\sO(\frac{d}{T}).
\end{align}
\end{theorem}

The proof of Theorem~\ref{theorem:main_convex} is provided in Appendix~\ref{appendix:proof_convex}. Some key implications  of Theorem~\ref{theorem:main_convex}  are discussed below:

\setdefaultleftmargin{12pt}{0pt}{}{}{}{}
\begin{enumerate}
\item  {\bf Convergence rate}: The dominant term on the right-hand side of~\eqref{eq:convergence_convex}  is $\sO(d/T^{\frac{2}{3}})$, which is better than  the convergence rate  $\sO(d^{{\frac{2}{3}}}/T^{\frac{1}{3}})$ of the vanilla SZO method  \cite{gasnikov2017stochastic} whenever $T>d$. Compared with the residual-feedback SZO method  with the rate of $\sO(d^2/T^{\frac{2}{3}})$ given in \cite{zhang2021new}, the  HLF-SZO \eqref{eq:single:both} has the same dependency on $T$ but achieves better dependency on the problem dimension $d$, which is due to the refined analysis on the second moment of the zeroth-order gradient estimator. On the other hand, all these SZO methods are inferior to the two-point ZO methods~\cite{nesterov2017random} that achieve the  convergence rate of $\sO(d/T)$. See Table \ref{tab:complexility} for the iteration complexity of these ZO methods.
\item {\bf Choice of $\beta$}: It is seen that the optimal $\beta$ that minimizes the dominant  term in~\eqref{eq:convergence_convex} is given by $\beta=1$. This choice of $\beta$ implied by the theoretical analysis is consistent with our empirical results in Section~\ref{sec:numerical}. As mentioned above, the HF-SZO method \eqref{eq:single:high} with $\beta=1$ is equivalent to the residual-feedback SZO method \cite{zhang2021new}.
\item {\bf Choice of $\alpha$}: It is seen that the parameter $\alpha$ does not appear in the dominant  term on the right-hand side of~\eqref{eq:convergence_convex}.
This suggests that, theoretically, the HLF-SZO \eqref{eq:single:both} converges at least as fast as the HF-SZO \eqref{eq:single:high} (or the residual-feedback SZO).
 Moreover, the numerical results in Section~\ref{sec:numerical} show that the HLF-SZO \eqref{eq:single:both} with $\alpha>0$ empirically achieves faster convergence than the residual-feedback SZO \cite{zhang2021new}.
\end{enumerate}

Next, we establish the convergence of HLF-SZO~\eqref{eq:single:both} for a nonconvex objective function $f$ with Theorem \ref{theorem:main_nonconvex}.

\begin{theorem}\label{theorem:main_nonconvex}
(Nonconvex Case). Let $\alpha\in[0,1)$,  $\beta\in(0,2)$, and the total number of iterations $T$ be sufficiently large. Suppose that Assumption \ref{ass:lip} holds and 
$f^\ast\coloneqq \inf_{x\in\mathbb{R}^d} f(x)>-\infty$. Then by choosing \eqref{eq:condition_convergence},
the HLF-SZO method \eqref{eq:single:both} achieves
\begin{align} \label{eq:con_nonconvex_pre}
 &   \frac{1}{T}\sum_{k=1}^T\mathbb{E}\!\left[\|\nabla f(\bx_k)\|^2\right] \nonumber \\
& \ \ \leq
\frac{4(1\!-\!\alpha)(f(\bx_1)-f^\ast)}{\eta T}
+\frac{8 G^2}{T^{2/3}}
+\sO(\frac{d}{T}).
\end{align}
Moreover, by letting $\eta$ achieve the equality in~\eqref{eq:condition_convergence}, we have
\begin{align}\label{eq:convergence_nonconvex}
& \ \ \frac{1}{T}\sum_{k=1}^T\mathbb{E}\!\left[\|\nabla f(\bx_k)\|^2\right]\nonumber\\
& \leq
\frac{d}{T^{2/3}}
\left(
\frac{80 L(f(\bx_1)-f^\ast)}{(1\!-\!|1\!-\!\beta|)^2}
+\frac{8 G^2}{d}
\right)
+\sO(\frac{d}{T}).
\end{align}
\end{theorem}
The proof of Theorem~\ref{theorem:main_nonconvex} is provided in Appendix~\ref{appendix:proof_nonconvex}. Note that in \eqref{eq:con_nonconvex_pre} and \eqref{eq:convergence_nonconvex}, we use the ergodic rate of the squared norm of the gradient $\|\nabla f(\bx_k)\|^2$ to characterize the convergence behavior. This is common for local methods of unconstrained smooth nonconvex problems, where one does not aim for global optimal solutions \cite{ghadimi2013stochastic,nesterov2017random}. Besides, it is observed that the bound~\eqref{eq:convergence_nonconvex} is very similar to~\eqref{eq:convergence_convex}, and   the discussions on Theorem~\ref{theorem:main_convex} can also be applied to the smooth nonconvex case with minor modifications.

\section{Numerical Experiments} \label{sec:numerical}

In this section, we first  study the properties of the  proposed SZO methods that incorporate high-pass and low-pass filters. 
Then   the   HLF-SZO \eqref{eq:single:both} is tested on logistic regression, ridge regression, an  artificial test function, and the linear quadratic regulator (LQR) problem,  in the comparison with the  residual-feedback SZO and the two-point ZO method.

\subsection{Properties of SZO Methods with  Filters}\label{sec:sim:prop}

This subsection compares the performance of  vanilla SZO \eqref{eq:single:basic}, HF-SZO \eqref{eq:single:high}, LF-SZO \eqref{eq:single:low}, and HLF-SZO \eqref{eq:single:both}, and studies the impact of the parameter $\beta$ on HLF-SZO \eqref{eq:single:both},
via the numerical tests on logistic regression. 


 Consider solving the logistic regression problem \eqref{eq:logis} \cite{uribe2020dual}:
\begin{align} \label{eq:logis}
    \min_{\bx\in\R^d} f(\bx) = \frac{1}{N} \sum_{i=1}^N \log\big(  1+\exp(-y_i\cdot A_i^\top \bx) \big),
\end{align}
where  $A_i\in \R^d$ is one of the data samples, $y_i\in\{-1,1\}$ is the corresponding class label, and $N$ is the total sample size. In our experiments, each element of a data sample $A_i$ is  randomly and independently generated from the uniform distribution $\mathrm{Unif}([-1,1])$, and the label is computed by $y_i=\mathrm{sign}(A_i^\top \bx_\star+\epsilon_i )$ with $\bx_\star = 0.5\mathbf{1}_d$ and $\epsilon_i\sim \mathrm{Unif}([-0.5,0.5])$.

\begin{figure}
    \centering
  \includegraphics[scale = 0.183]{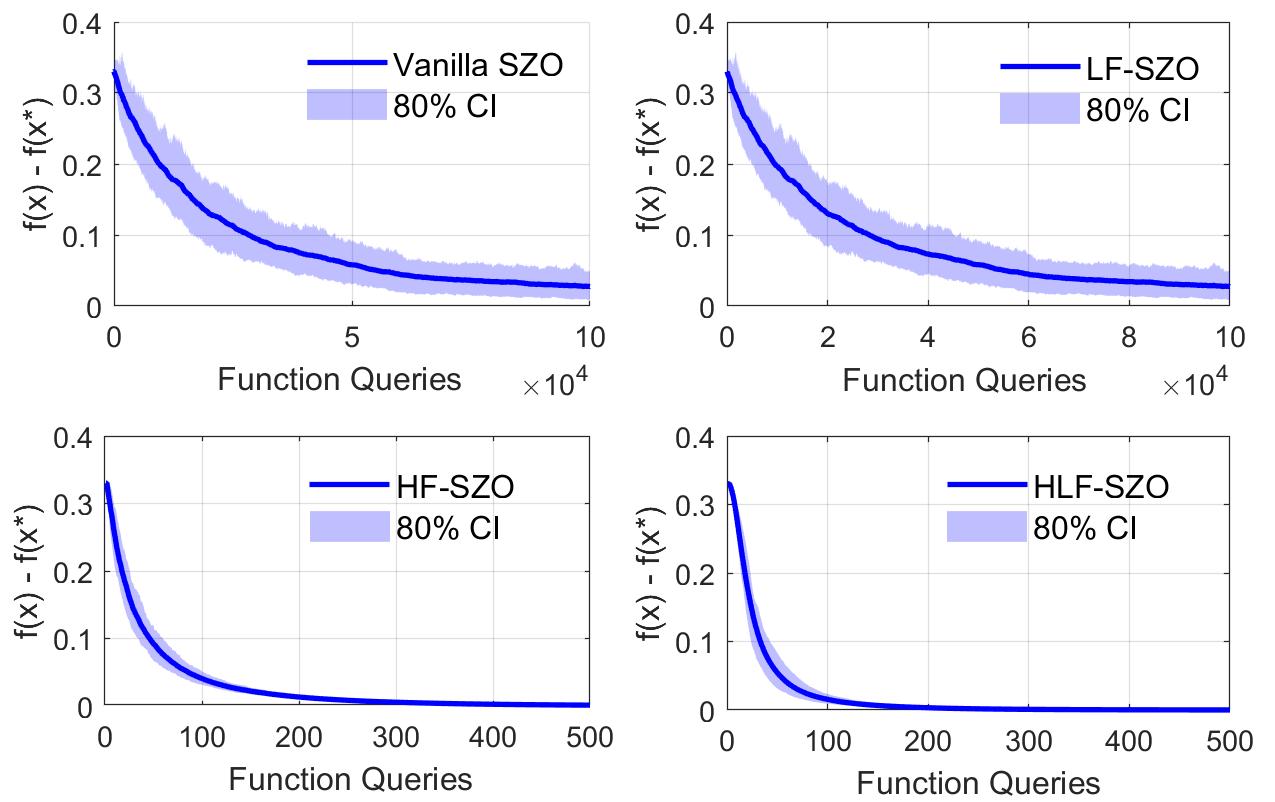}
    \caption{The convergence results of applying vanilla SZO \eqref{eq:single:basic}, HF-SZO \eqref{eq:single:high}, LF-SZO \eqref{eq:single:low}, and HLF-SZO \eqref{eq:single:both} to solve the logistic regression \eqref{eq:logis} in Case 1a).}
    \label{fig:logiscom}
\end{figure}
\begin{figure}
    \centering
  \includegraphics[scale = 0.225]{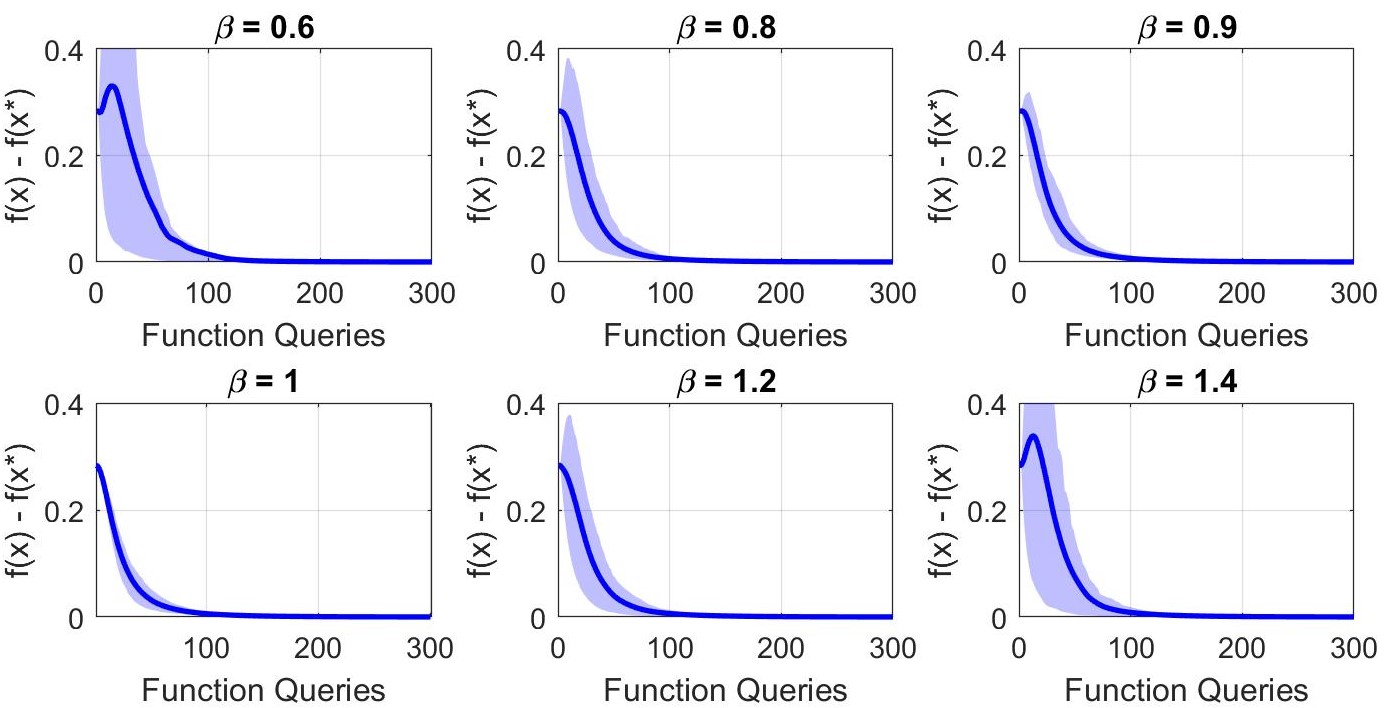}   
    \caption{The convergence results of  HLF-SZO \eqref{eq:single:both} under different $\beta$ for solving the logistic regression \eqref{eq:logis} in Case 1b).}
    \label{fig:log:beta}
\end{figure}

\textbf{Case 1a)}. Let $d=2$ and $N=200$. We set $r=0.1$, $\beta = 1$, $\alpha=0.9$ and the initial condition $\bx_0=\bm{0}_d$ for the SZO methods. 
Then we  manually optimize the stepsize $\eta$ of each SZO method to achieve its fastest convergence\footnote{We gradually tune the stepsize for each SZO method and select the one that has the fastest convergence without divergence.}. The selected stepsizes are $5\times 10^{-4}$,  $0.3$, $5\times 10^{-5}$, and $0.05$ for the vanilla SZO \eqref{eq:single:basic}, HF-SZO \eqref{eq:single:high}, LF-SZO \eqref{eq:single:low},  and HLF-SZO
 \eqref{eq:single:both}, respectively.  We run each SZO method 200 times and calculate the mean and $80\%$-confidence interval (CI) of the distance to optimality, i.e., $f(\bx_k)-f(\bx^*)$.
 Figure \ref{fig:logiscom} illustrates the experimental results of solving the logistic regression \eqref{eq:logis} with these SZO methods.

\textbf{Case 1b)}. 
Next, we 
 tune the parameter $\beta$ from 0.6 to 1.4 and run  experiments for each case to study the impact of $\beta$. Other settings are the same as the above Case 1a).
 The convergence results of HLF-SZO \eqref{eq:single:both} with different $\beta$ are  shown in Figure \ref{fig:log:beta}.

The key observations are summarized as follows:
\begin{itemize}
    \item [1)] From Figure \ref{fig:logiscom}, it is seen that  HF-SZO and HLF-SZO have much  smaller
     variances and converge much faster than the vanilla SZO. It indicates that the integration of a high-pass filter can greatly reduce the estimation variance of SZO, and thus leads to faster convergence.
 
  \item [2)]  Figure \ref{fig:logiscom} also shows  that HLF-SZO  converges even faster  than HF-SZO due to the integration of a low-pass filter. It implies that the momentum term introduced by the low-pass filter can further accelerate convergence, which is verified by  more numerical tests in
  the next subsections.
    
    \item [3)] From Figure \ref{fig:log:beta}, it is observed that the case with  $\beta =1$ achieves the least variance and fastest convergence. This is consistent with the theoretical analysis in Section \ref{sec:analysis} that $\beta=1$ is the optimal choice.
    
\end{itemize}

\subsection{Comparison with  Other ZO Algorithms}

This subsection compares the performance of the proposed HLF-SZO method \eqref{eq:single:both} with the  residual-feedback SZO method  \cite{zhang2021new} and the two-point ZO method \eqref{eq:two}, through the tests on logistic regression, ridge regression, and Beale function.
Here,
the vanilla SZO method  \eqref{eq:single:basic} is not considered for comparison, because it has a much larger variance and much slower convergence.

\begin{figure}
   \centering
    \includegraphics[scale=0.345]{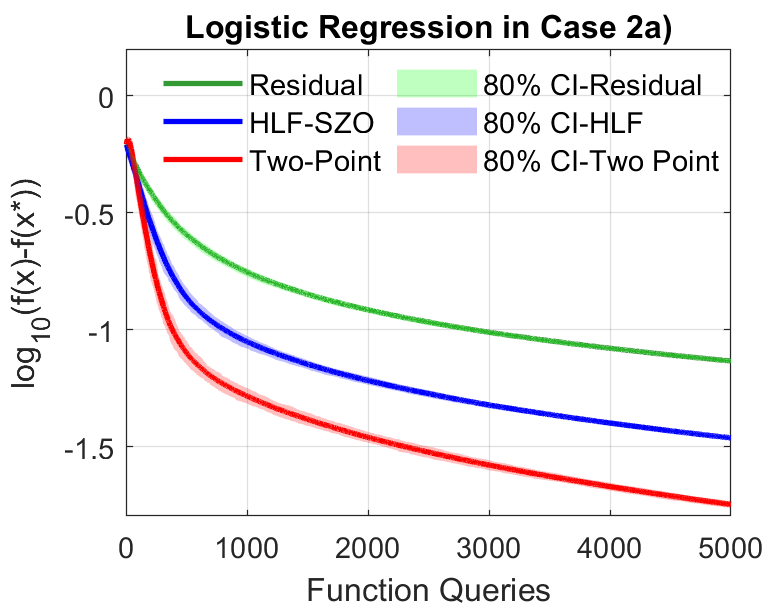}
      \includegraphics[scale=0.35]{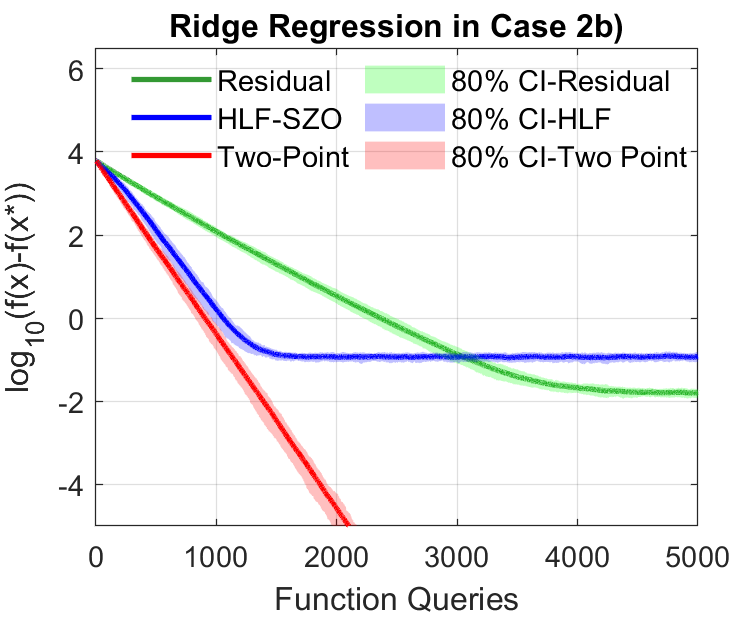}
        \includegraphics[scale=0.345]{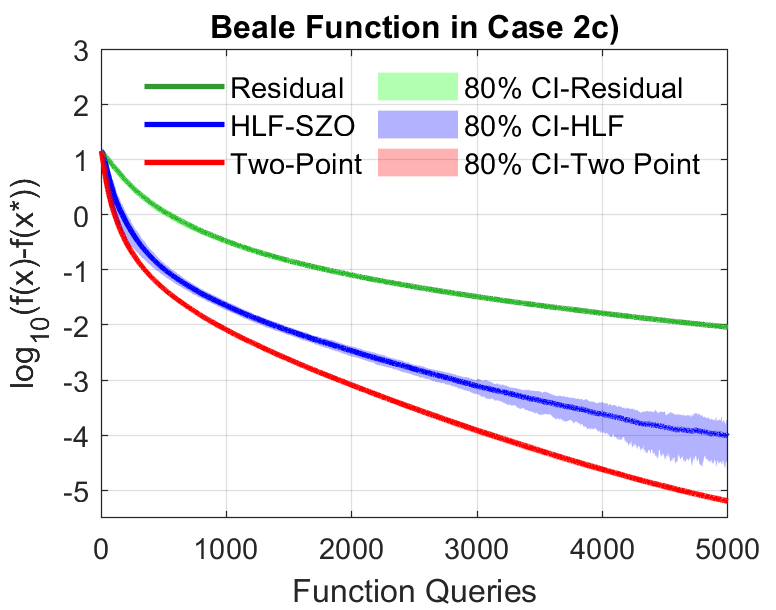}
    \caption{The convergence results of residual-feedback SZO, HLF-SZO \eqref{eq:single:both}, and two-point ZO \eqref{eq:two} for solving  {logistic regression} \eqref{eq:logis} in Case 2a), ridge regression \eqref{eq:ridge} in Case 2b), and Beale function \eqref{eq:beale} in Case 2c).
    }
    \label{fig:convergence:comp}
\end{figure}

In the following cases,
we run each ZO method 200 times and calculate the mean and $80\%$-confidence interval (CI) of the distance to optimality, i.e., $f(\bx_k)-f(\bx^*)$. Figure~\ref{fig:convergence:comp} illustrates  the convergence results of these ZO methods.




\textbf{Case 2a)}. Consider solving the  \textbf{logistic regression} problem \eqref{eq:logis} with   $d=50$ and $N=1000$. We manually optimize the stepsize $\eta$ of each ZO method to achieve its fastest convergence. The selected  stepsizes are $1.5\times 10^{-2}$,   $4.5\times 10^{-2}$,  and $0.7$ for the HLF-SZO \eqref{eq:single:both}, the residual-feedback SZO, and the two-point ZO \eqref{eq:two}.
Other settings are the same as Case 1a) in Section \ref{sec:sim:prop}. 



\textbf{Case 2b)}.  Consider solving the   \textbf{ridge regression} problem  \eqref{eq:ridge} \cite{uribe2020dual}:
\begin{align}\label{eq:ridge}
    \min_{\bx\in\R^d}\, f(\bx)= \frac{1}{2}||\bm{b} - H\bx ||_2^2 + \frac{1}{2}c ||\bx||_2^2.
\end{align}
In our experiments,
each entry of matrix $H\in \R^{N\times d}$ is generated randomly and independently from  Gaussian distribution $\mathcal{N}(0,1)$. The vector $\bm{b}\in \R^N$ is constructed by letting
$\bm{b}= H\bx_{\star} +\bm{\epsilon}$ for a certain predefined   $\bx_\star\in \R^d$ and $\bm{\epsilon}\sim \mathcal{N}(\bm{0},0.1\bm{I})$.  The regularization constant is set to be $c=0.1$. Let $d=50$, $N=1000$ and $\bx_\star = 0.5\times \bm{1}_d$. For the ZO methods, we set $r=0.1$,  $\alpha=0.9$, $\beta =1$, and manually optimize the stepsize $\eta$ to achieve the fastest convergence. The selected stepsizes are  $1\times 10^{-6}$,   $2.4\times 10^{-6}$, and $2\times 10^{-5}$  for the HLF-SZO \eqref{eq:single:both}, the residual-feedback SZO, and the two-point ZO \eqref{eq:two}, respectively. 


\textbf{Case 2c)}. Consider minimizing the  nonconvex \textbf{Beale function} 
 given by \eqref{eq:beale}, which is an  artificial   test function and
has a global minimum at $\bx^*=[3;0.5]$ with $f(\bx^*)=0$. 
\begin{align}
     f(\bx) &= (1.5-x_1+x_1x_2)^2   + (2.25\!-\!x_1+x_1x_2^2)^2 \nonumber\\
     &\qquad\qquad  + (2.625\!-\!x_1+x_1x_2^3)^2.  \label{eq:beale}
\end{align}
For the ZO methods, we set $r=0.01$, $\beta=1$, $\alpha=0.9$, $\bx_0=[0;0]$, and manually optimize the stepsize $\eta$ to achieve the fastest convergence. The selected stepsizes are $2\times 10^{-4}$,   $5.8\times 10^{-4}$, and $6\times 10^{-3}$  for the HLF-SZO \eqref{eq:single:both}, the residual-feedback SZO, and the two-point ZO \eqref{eq:two},  respectively. 



From Figure \ref{fig:convergence:comp}, it is seen that the proposed HLF-SZO \eqref{eq:single:both} converges faster 
 than the residual-feedback SZO method in all cases, and its performance is comparable to that of  the two-point ZO method \eqref{eq:two}. 
Since the residual-feedback SZO  is equivalent to the HF-SZO \eqref{eq:single:high} with $\beta =1$, 
these results indicate that the integration of a low-pass filter can further accelerate the convergence of SZO via the momentum term. 

\subsection{Control Policy Optimization for Linear Quadratic Regulator (LQR) Problem}

We consider solving the  linear quadratic regulator (LQR) problem \cite{anderson2007optimal,fazel2018global} \eqref{eq:lqr} with the linear feedback control policy $\bu_t = - K\bx_t$.
\begin{subequations}\label{eq:lqr}
\begin{align}
    \min_K \quad & J(K)\coloneqq \mathbb{E}\Big[\sum_{t=0}^\infty \gamma^t\big( \bx_t^\top Q\bx_t + \bu_t^\top R \bu_t  \big)\Big]\\
    \text{s.t.}\quad  & \bx_{t+1} = A \bx_t +B \bu_t + \bw_t.
\end{align}
\end{subequations}
Here, $\bx_t\in\R^{n_x}$ and $\bu_t\in\R^{n_u}$ are the  state and the control input at time step $t$, respectively. $\bw_t\sim \mathcal{N}(0, \delta^2 I)$ denotes the random system noise.  
In the simulations, we set $n_x = 20$, $n_u=15$,  and thus the problem dimension is $d= 300$. Let the discount factor $\gamma$ be 0.9.
Each element in the transition  matrices $A\in \R^{n_x\times n_x}$ and $B\in\R^{n_x\times n_u}$ are randomly generated
from $\text{Unif}([-0.1,0])$ and $\text{Unif}([0,0.02])$, respectively.
Let 
  $Q$ and $R$ be the identity matrix.
To evaluate the total cost $J(K)$, we run the episode with a finite horizon $H=50$ to approximate it, which has been checked to be sufficiently long  for the system to converge. The initial state is set as $\bx_0=\bm{1}_{n_x}$.

We apply the  HLF-SZO \eqref{eq:single:both},  residual-feedback SZO, and  two-point ZO \eqref{eq:two} to solve
the LQR problem \eqref{eq:lqr}.  The initial $K_0$ is constructed by $K_0 = K^* + \tilde{K}$, where $K^*$ is the  optimal solution and each element in $\tilde{K}$ is randomly  generated from $\text{Unif}([-1,1])$.
We  set $r=0.1$,  $\alpha=0.9$, $\beta =1$, and tune the standard deviation $\delta$ of $\bw_t$ to simulate different levels of system noise, as in the following cases. In each case, we
manually optimize the stepsize $\eta$ of each ZO method to achieve its fastest convergence.

\begin{figure}
   \centering
      \includegraphics[scale = 0.345]{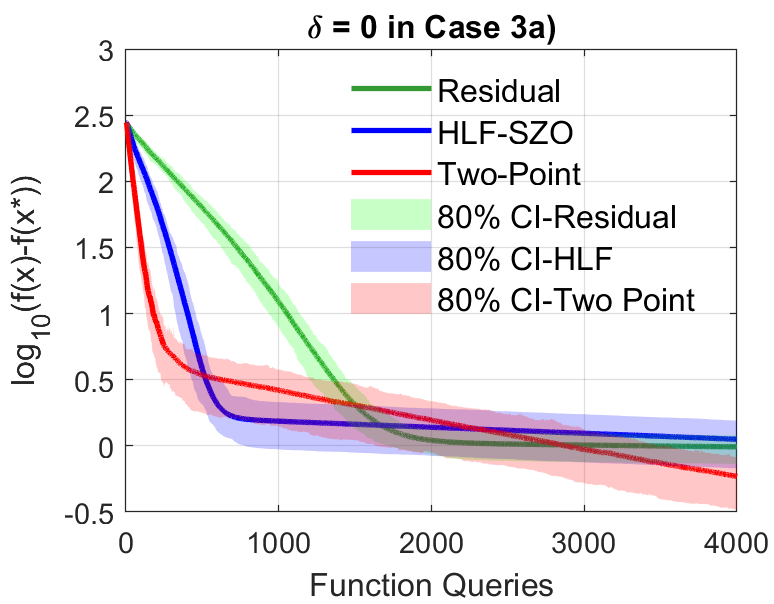}   
    \includegraphics[scale = 0.345]{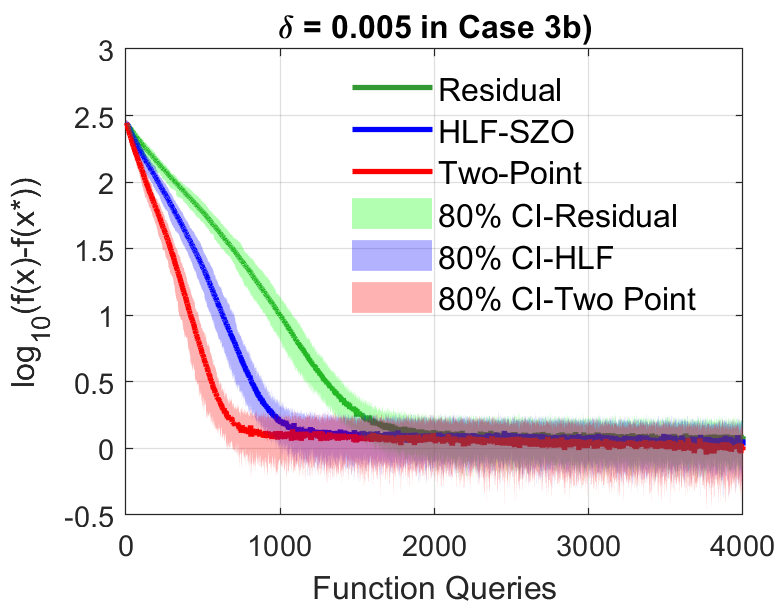}
     \includegraphics[scale=0.347]{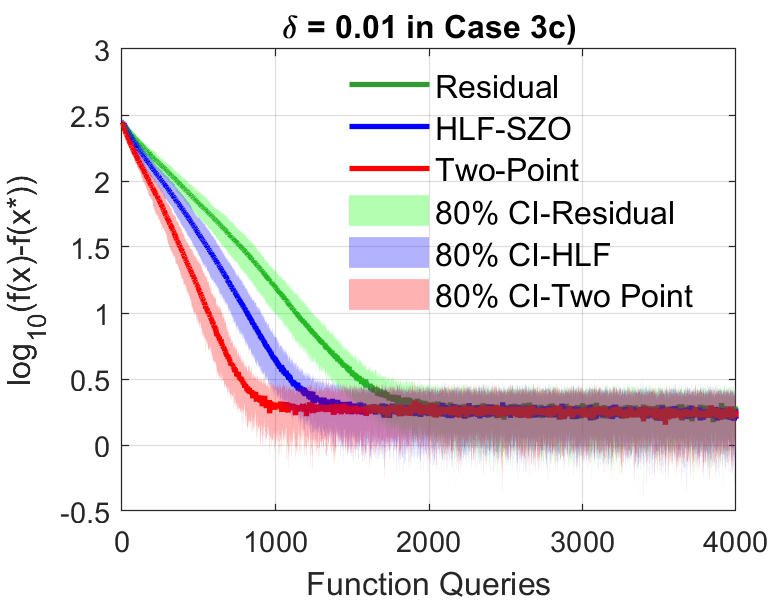}
    \caption{The convergence results of residual-feedback SZO, HLF-SZO \eqref{eq:single:both}, and two-point ZO \eqref{eq:two} for solving  the LQR problem \eqref{eq:lqr} under different levels of system noise $\delta$.
    }
    \label{fig:lqr}
\end{figure}

\textbf{Case 3a)}. Let  $\delta=0$.
The selected stepsizes are  $7.5\times 10^{-6}$,   $2.25\times 10^{-5}$, and $7.5\times 10^{-4}$  for  HLF-SZO \eqref{eq:single:both},  residual-feedback SZO, and  two-point ZO \eqref{eq:two}, respectively.   

\textbf{Case 3b)}. Let  $\delta=0.005$. The selected stepsizes are  $4\times 10^{-6}$,   $2.4\times 10^{-5}$, and $1.2\times 10^{-4}$  for  HLF-SZO \eqref{eq:single:both},  residual-feedback SZO, and  two-point ZO \eqref{eq:two}, respectively.

\textbf{Case 3c)}. Let  $\delta=0.01$. The selected stepsizes are  $3\times 10^{-6}$,   $2.1\times 10^{-5}$, and $9\times 10^{-5}$  for  HLF-SZO \eqref{eq:single:both},  residual-feedback SZO, and  two-point ZO \eqref{eq:two}, respectively.

The convergence results of these ZO methods are illustrated in Figure \ref{fig:lqr}. It is observed that larger noise leads to slower convergence  and lower convergence accuracy of ZO methods.
In all the  cases, the HLF-SZO \eqref{eq:single:both} can converge faster than the residual-feedback SZO and achieve comparable  performance to the two-point ZO method \eqref{eq:two}.

\section{Conclusion}

In this work, by leveraging the idea of high-pass and low-pass filters from extremum seeking control, we develop a  novel single-point zeroth-order optimization (SZO) method called HLF-SZO. Interestingly, we find that the integration of a high-pass filter coincides with the residual feedback scheme proposed in \citet{zhang2021new}, and  the integration of a low-pass filter can be interpreted as the momentum  method. 
We prove that the  HLF-SZO method  achieves an iteration complexity of  $\sO(d^{\frac{3}{2}}/\epsilon^{\frac{3}{2}})$ for both convex and nonconvex objective functions under the Lipschitz and smoothness conditions. This iteration complexity is improved over the vanilla SZO method   \cite{gasnikov2017stochastic} and the residual-feedback SZO method \cite{zhang2021new}, although it is inferior to  the   two-point ZO methods. Extensive 
numerical experiments demonstrate that the HLF-SZO method  achieves  a much smaller variance and much faster convergence than the vanilla SZO method; it  empirically outperforms the residual-feedback SZO method and  has  comparable  performance to  the two-point ZO method.

Bridging ZO  and  extremum seeking control  
to benefit each other is a promising direction of research. This paper is a preliminary attempt to throw light on this direction, and much remains to be studied. Our future work includes  
1) extending HLF-SZO to stochastic optimization problems, 2) designing better  filters and time-discretization schemes to  further improve convergence, etc.

\section*{Acknowledgements}

This work was supported by the funding programs:    NSF CAREER ECCS-1553407, 
NSF AI Institute 2112085, NSF CNS 2003111, and
ONR YIP N00014-19-1-2217.


\bibliography{Reference}
\bibliographystyle{icml2022}

\newpage
\appendix
\onecolumn
\section{Notations and Preliminary Results}

It can be checked that the iterations can be equivalently written as
\begin{align*}
\bg_k &= \frac{d}{r}z_k \bu_k,\quad
z_k = (1-\beta) z_{k-1}
+ f(\bx_k\!+\!r\bu_k) - f(\bx_{k-1}\!+\!r\bu_{k-1}), \\
\bp_k &= \alpha\bp_{k-1}+\eta\bg_k, \\
\bx_{k+1} &= \bx_k - \bp_k,
\end{align*}
for $k\geq1$, and we set $\bp_0=\bm{0}$, $\bx_1=\bx_0$ and $z_0=0$. We denote $\tilde{\beta}\coloneqq 1-|1-\beta|$ for notational simplicity. We treat $\bx_1=\bx_0$ as deterministic quantities and let $\mathcal{F}_k$ denote the filtration generated by $\bx_2,\ldots,\bx_k$ and $\bu_0,\bu_1,\ldots,\bu_{k-1}$ for $k\geq 1$.

The following lemma deals with the expectation of $\bg_k$.
\begin{lemma}\label{lemma:function_fr}
Denote $\mathbb{B}_d\coloneqq \{\bx\in\mathbb{R}^d:\|\bx\|\leq 1\}$. For all $k\geq 1$, we have 
$$
\mathbb{E}[\bg_k\,|\,\mathcal{F}_k]
=\nabla f_r(\bx_k),
$$
where $f_r:\mathbb{R}^d\rightarrow\mathbb{R}$ is defined by
$
f_r(\bx) \coloneqq \mathbb{E}_{\by\sim\mathrm{Unif}(\mathbb{B}_d)}[f(\bx+r\by)]
$. Moreover, for all $\bx\in\mathbb{R}^d$,
$$
|f_r(\bx)-f(\bx)|\leq\frac{1}{2}Lr^2,
\qquad
\|\nabla f_r(\bx)-\nabla f(\bx)\|
\leq Lr,
$$
and if $f$ is convex, then $f_r(\bx)\geq f(\bx)$.
\end{lemma}
\begin{proof}
The result $\mathbb{E}[\bg_k\,|\,\mathcal{F}_k]=\nabla f_r(\bx_k)$ is standard in zeroth-order optimization~\cite{flaxman2005online}. By the $L$-smoothness of $f$, we have
$$
r\langle \by,\nabla f(\bx)\rangle
-\frac{L}{2}r^2\|\by\|^2\leq f(\bx+r\by) - f(\bx)
\leq r\langle \by,\nabla f(\bx)\rangle
+\frac{L}{2}r^2\|\by\|^2,
$$
and by taking the expectation with $\by\sim\mathrm{Unif}(\mathbb{B}_d)$, we get
$$
-\frac{L}{2}r^2\leq f_r(\bx) - f(\bx)
\leq \frac{L}{2}r^2
$$
for any $\bx\in\mathbb{R}^d$. If $f$ is convex, then $f(\bx+r\by)-f(\bx)\geq r\langle \by,\nabla f(\bx)\rangle$, and by by taking the expectation with $\by\sim\mathrm{Unif}(\mathbb{B}_d)$, we get $f_r(\bx)\geq f(\bx)$.

Finally, regarding the gradient $\nabla f_r(\bx)$, we have
\begin{align*}
\|\nabla f_r(\bx) - \nabla f(\bx)\|
=\ &
\|\nabla \mathbb{E}_{\by\sim\mathrm{Unif}(\mathbb{B}_d)}[f(\bx+r\by)] - \nabla f(\bx)\| \\
=\ &
\|\mathbb{E}_{\by\sim\mathrm{Unif}(\mathbb{B}_d)}[\nabla f(\bx+r\by) - \nabla f(\bx)]\| \\
\leq\ &
\mathbb{E}_{\by\sim\mathrm{Unif}(\mathbb{B}_d)}
[Lr\|\by\|]\leq Lr,
\end{align*}
where we interchange integration and differentiation in the second step, and it can be justified by the dominated convergence theorem.
\end{proof}

The following lemma deals with the accumulated second moment of $z_k$.
\begin{lemma}\label{lemma:zk_second_moment}
We have
$$
\sum_{k=1}^T\mathbb{E}\!\left[
|z_k|^2\right]
\leq
\frac{5r^2}{\tilde{\beta}^2d}
\sum_{k=1}^{T-1}
\mathbb{E}\!\left[
\|\nabla f(\bx_{k})\|^2
\right]
+\frac{2G^2}{\tilde{\beta}^2}
\sum_{k=1}^{T-1}
\mathbb{E}\!\left[\|\bp_{k}\|^2\right]
+\frac{10Tr^4L^2}{\tilde{\beta}^2}
+\frac{5r^2}{2\tilde{\beta}d}\|\nabla f(\bx_1)\|^2
.
$$
\end{lemma}
\begin{proof}
For $k=1$, by the initial conditions imposed on $\bx_0$ and $z_0$, we have
\begin{align*}
\mathbb{E}\!\left[|z_1|^2\right]
=\ &
\mathbb{E}\!\left[|f(\bx_1\!+\!r\bu_1) - f(\bx_1 \!+\! r\bu_0)|^2\right] =   \mathbb{E}\!\left[
\left|f(\bx_1\!+\!r\bu_1) -  f(\bx_1) -  ( f(\bx_1 \!+\! r\bu_0)-f(\bx_1)) \right|^2\right]         \\
=\ &    \mathbb{E}\!\left[
\left| \int_{0}^r\langle \nabla f(\bx_1 \!+\! s\bu_1),\bu_1\rangle\,ds -  \int_{0}^r\langle \nabla f(\bx_1 \!+\! s\bu_0),\bu_0\rangle\,ds  \right|^2\right]         \\
=\ &
\mathbb{E}\!\left[
\left|
r\langle \nabla f(\bx_1), \bu_1 \!-\! \bu_0\rangle
+\int_{0}^r\langle \nabla f(\bx_1 \!+\! s\bu_1) \!-\! \nabla f(\bx_1),\bu_1\rangle\,ds
-\int_{0}^r\langle \nabla f(\bx_1 \!+\! s\bu_0) \!-\! \nabla f(\bx_1),\bu_0\rangle\,ds
\right|^2\right] \\
\leq\ &
\frac{5}{4}r^2\nabla f(\bx_1)^{\!\top}\mathbb{E}\!\left[(\bu_1-\bu_0)(\bu_1-\bu_0)^\top\right]\nabla f(\bx_1) \\
&
+10\,
\mathbb{E}\!\left[
\left|\int_{0}^r
\|\nabla f(\bx_1 \!+\! s\bu_1)-\nabla f(\bx_1)\| \|\bu_1\|\,ds
\right|^2
\right]
+10\,
\mathbb{E}\!\left[
\left|\int_{0}^r
\|\nabla f(\bx_1 \!+\! s\bu_0)-\nabla f(\bx_1)\| \|\bu_0\|\,ds
\right|^2
\right]\\
\leq\ &
\frac{5r^2}{2d}\|\nabla f(\bx_1)\|^2
+10
\mathbb{E}\!\left[
\left|
\int_{0}^r L|s|\|\bu_1\|^2\,ds\right|^2
\right]
+10
\mathbb{E}\!\left[
\left|
\int_{0}^r L|s|\|\bu_0\|^2\,ds\right|^2
\right] \\
\leq\ &
\frac{5r^2}{2d}\|\nabla f(\bx_1)\|^2
+5L^2 r^4.
\end{align*}
The first inequality above is the application of the inequality $(a+b+c)^2 \leq \frac{5}{4}a^2 + 10b^2 +10c^2$, as $\frac{5}{4}a^2 + 10b^2 +10c^2 - (a+b+c)^2 = (\frac{1}{\sqrt{8}}a - \sqrt{8}b)^2 + (\frac{1}{\sqrt{8}}a - \sqrt{8}c)^2 + (b-c)^2\geq 0$. The second inequality above is due to $\mathbb{E}\!\left[(\bu_1-\bu_0)(\bu_1-\bu_0)^\top\right] = \frac{2}{d}I_d$.

Then for $k>1$,  by the inequality  $(a+b)^2\leq  (1+\gamma) a^2 +(1+\frac{1}{\gamma})b^2 $ for any $\gamma>0$ (frequently used below to bound the square of a sum),  we let $\gamma = \frac{1-|1-\beta|}{|1-\beta|}$ and get
\begin{align*}
\mathbb{E}\!\left[|z_k|^2\right]
\leq\ &
\frac{1}{|1\!-\!\beta|}\cdot(1\!-\!\beta)^2\,
\mathbb{E}\!\left[|z_{k-1}|^2\right]
+\frac{1}{1\!-\!|1\!-\!\beta|}\,
\mathbb{E}\!\left[\left|
f(\bx_k\!+\!r\bu_k)
\!-\!
f(\bx_{k-1}\!+\!r\bu_{k-1})\right|^2\right] \\
=\ &
(1\!-\!\tilde{\beta})\,
\mathbb{E}\!\left[|z_{k-1}|^2\right]
+
\tilde{\beta}^{-1}
\mathbb{E}\!\left[\left|
f(\bx_k\!+\!r\bu_k)
\!-\!
f(\bx_{k-1}\!+\!r\bu_{k-1})\right|^2\right]
\end{align*}
when $\beta\neq 1$. And obviously this inequality also extends to the $\beta=1$ case, because $z_k = f(\bx_k\!+\!r\bu_k)
\!-\!
f(\bx_{k-1}\!+\!r\bu_{k-1})$ and $\tilde{\beta}=1$ when $\beta =1$. Furthermore,
\begin{align*}
&\mathbb{E}\!\left[
\left|f(\bx_k \!+\! r\bu_k) \!-\! f(\bx_{k-1} \!-\! r\bu_{k-1})\right|^2
\right] \\
=&  \mathbb{E}\!\left[
\left|f(\bx_k \!+\! r\bu_k) \!- f(\bx_{k-1}\!+\!r\bu_{k}) + f(\bx_{k-1}\!+\!r\bu_{k})-\! f(\bx_{k-1} \!-\! r\bu_{k-1})\right|^2
\right]  \\
\leq\   
&
2
\Big(
\mathbb{E}\!\left[
\left|f(\bx_{k-1}\!+\!r\bu_k) \!-\! f(\bx_{k-1}\!+\!r\bu_{k-1})\right|^2\right]   +
\mathbb{E}\!\left[
\left|f(\bx_{k}\!+\!r\bu_{k}) \!-\! f(\bx_{k-1}\!+\!r\bu_{k})\right|^2\right]\!
\Big).
\end{align*}
For the first term, we have
\begin{align*}
& \mathbb{E}\!\left[|f(\bx_{k-1}\!+\!r\bu_k)-f(\bx_{k-1}+r\bu_{k-1})|^2\right] \\
=\ &
\mathbb{E}\Bigg[
\bigg|
r\langle \nabla f(\bx_{k-1}), \bu_k \!-\! \bu_{k-1}\rangle
+\int_0^r\langle \nabla f(\bx_{k-1}\!+\!s\bu_k)\!-\!\nabla f(\bx_{k-1}),\bu_k\rangle\,ds \\
&
\quad\ \ -\int_0^r\langle \nabla f(\bx_{k-1}\!+\!s\bu_{k-1})\!-\!\nabla f(\bx_{k-1}),\bu_{k-1}\rangle\,ds
\bigg|^2\Bigg] \\
\leq\ &
5\,\mathbb{E}\!\left[
\left(
\int_0^r\! \left(
\|\nabla f(\bx_{k-1}\!+\!s\bu_k)\!-\!\nabla f(\bx_{k-1})\| + \|\nabla f(\bx_{k-1}\!+\!s\bu_{k-1})\!-\!\nabla f(\bx_{k-1})\|\right) ds\right)^2
\right] \\
&
+\frac{5}{4}r^2\,\mathbb{E}\!\left[|\langle\nabla f(\bx_{k-1}),\bu_k\!-\!\bu_{k-1}\rangle|^2\right] \\
\leq\ &
5\,
\mathbb{E}\left[
\left(\int_0^r
(Ls\|\bu_k\| + Ls \|\bu_{k-1}\|)\,ds\right)^2\right] \\ & +\frac{5}{4}r^2\,
\mathbb{E}\!\left[
\nabla f(\bx_{k-1})^{\!\top}
\mathbb{E}\!\left[\left.
(\bu_k\!-\!\bu_{k-1})(\bu_k\!-\!\bu_{k-1})^{\!\top}\,\right|\mathcal{F}_{k-1}\right]
\nabla f(\bx_{k-1})
\right] \\
=\ &
5 r^4 L^2
+\frac{5r^2}{2d} \mathbb{E}\!\left[\|\nabla f(\bx_{k-1})\|^2\right]
\end{align*}
where we used $\|\bu_k\|=\|\bu_{k-1}\|=1$, the law of total expectation, and
$$
\mathbb{E}\!\left[\left.
(\bu_k\!-\!\bu_{k-1})(\bu_k\!-\!\bu_{k-1})^{\!\top}\,\right|\mathcal{F}_{k-1}\right] = \mathbb{E}\!\left[\bu_k \bu_k^\top + \bu_{k-1} \bu_{k-1}^\top |\mathcal{F}_{k-1}\right]
=\frac{2}{d}I_d,
$$
and for the second term, we have
$$
\mathbb{E}\!\left[|f(\bx_k\!+\!r\bu_k)
-f(\bx_{k-1}\!+\!r\bu_k)|^2\right]
\leq G^2\,\mathbb{E}\left[\|\bx_k-\bx_{k-1}\|^2\right]
= G^2\,\mathbb{E}\left[\|\bp_{k-1}\|^2\right].
$$
Therefore for $k>1$,
\begin{align*}
\mathbb{E}\!\left[|z_k|^2\right]
\leq\ &
(1\!-\!\tilde{\beta})\,\mathbb{E}\!\left[|z_{k-1}|^2\right]
+\frac{5r^2}{\tilde{\beta} d}
\mathbb{E}\!\left[
\|\nabla f(\bx_{k-1})\|^2
\right]
+\frac{10r^4L^2}{\tilde{\beta}}
+\frac{2G^2}{\tilde{\beta}}
\mathbb{E}\!\left[\|\bp_{k-1}\|^2\right].
\end{align*}
By summing over $k=1,\ldots,T$, we get
\begin{align*}
\sum_{k=1}^T\mathbb{E}\!\left[|z_k|^2\right]\leq\ &
(1\!-\!\tilde{\beta})\sum_{k=1}^{T-1}\mathbb{E}\!\left[|z_{k}|^2\right]
+\frac{5r^2}{\tilde{\beta} d}
\sum_{k=1}^{T-1}\mathbb{E}\!\left[
\|\nabla f(\bx_{k})\|^2\right]
+\frac{10(T-1)r^4L^2}{\tilde{\beta}} \\
&
+\frac{2G^2}{\tilde{\beta}}\sum_{k=1}^{T-1}
\mathbb{E}\!\left[\|\bp_{k}\|^2\right]
+\frac{5r^2}{2d}\mathbb{E}\!\left[
\|\nabla f(\bx_1)\|^2\right]
+5r^4L^2 \\
\leq\ &
(1\!-\!\tilde{\beta})\sum_{k=1}^{T}\mathbb{E}\!\left[|z_{k}|^2\right]
+\frac{5r^2}{\tilde{\beta} d}
\sum_{k=1}^{T-1}\mathbb{E}\!\left[
\|\nabla f(\bx_{k})\|^2\right]
+\frac{10Tr^4L^2}{\tilde{\beta}}
+\frac{2G^2}{\tilde{\beta}}\sum_{k=1}^{T-1}
\mathbb{E}\!\left[\|\bp_{k}\|^2\right] +\frac{5r^2}{2d}\|\nabla f(\bx_1)\|^2
,
\end{align*}
which then implies the desired result.
\end{proof}

The following lemma then provides bounds for the accumulated second moment of $\bg_k$ and $\bp_k$.
\begin{lemma}\label{lemma:bound_gk_pk}
Suppose the quantity
$$
\theta\coloneqq 1-\frac{4\eta^2 d^2G^2(1+\alpha^2)}{\tilde{\beta}^2 r^2(1-\alpha^2)^2}
$$
is positive. Then we have
\begin{align}
\sum_{k=1}^T\mathbb{E}\!\left[\|\bg_k\|^2\right]
\leq\ &
\frac{5d}{\theta\tilde{\beta}^2}\sum_{k=1}^{T}\mathbb{E}\!\left[\|\nabla f(\bx_k)\|^2\right]
+ \frac{10Tr^2L^2d^2}{\theta\tilde{\beta}^2}
+\frac{5d}{2\theta\tilde{\beta}}\|\nabla f(\bx_1)\|^2,
\label{eq:bound_gk}
\\
\sum_{k=1}^T\mathbb{E}\!\left[\|\bp_k\|^2\right]
\leq\ &
\frac{2\eta^2(1+\alpha^2)}{(1-\alpha^2)^2}
\sum_{k=1}^T \mathbb{E}\!\left[\|\bg_k\|^2\right].
\label{eq:bound_pk}
\end{align}
\end{lemma}
\begin{proof}
By the definition of $\bp_k$, we have
\begin{align*}
\mathbb{E}\!\left[\|\bp_k\|^2\right]
\leq\ &
\left(1\!+\!\frac{1\!-\!\alpha^2}{2\alpha^2}\right)\alpha^2\mathbb{E}\!\left[\|\bp_{k-1}\|^2\right]
+ \left(1\!+\!\frac{2\alpha^2}{1\!-\!\alpha^2}\right)\eta^2\mathbb{E}\!\left[\|\bg_k\|^2\right] \\
=\ &
\frac{1\!+\!\alpha^2}{2}\mathbb{E}\!\left[\|\bp_{k-1}\|^2\right]
+\frac{1\!+\!\alpha^2}{1\!-\!\alpha^2}\eta^2
\mathbb{E}\!\left[\|\bg_k\|^2\right].
\end{align*}
By summing over $k=1,\ldots,T$, we get
$$
\sum_{k=1}^T \mathbb{E}\!\left[
\|\bp_k\|^2\right]
\leq
\frac{1+\alpha^2}{2}
\sum_{k=1}^{T-1}\mathbb{E}\!\left[
\|\bp_{k}\|^2
\right]
+\frac{1+\alpha^2}{1-\alpha^2}\eta^2\sum_{k=1}^T
\mathbb{E}\!\left[\|\bg_k\|^2\right],
$$
which then implies~\eqref{eq:bound_pk}.

Now, since $\bg_k=\frac{d}{r}z_k\bu_k$ and $\|\bu_k\|=1$, we see that
$\mathbb{E}\!\left[\|\bg_k\|^2\right]=\frac{d^2}{r^2}\mathbb{E}\!\left[|z_k|^2\right]$. By using the bound in Lemma~\ref{lemma:zk_second_moment}, we get
\begin{align*}
\sum_{k=1}^T \mathbb{E}\!\left[\|\bg_k\|^2\right]
\leq\ &
\sum_{k=1}^{T-1}\frac{2d^2 G^2}{\tilde{\beta}^2 r^2}\mathbb{E}\!\left[\|\bp_{k}\|^2\right]
+\frac{5d}{\tilde{\beta}^2}\sum_{k=1}^{T-1}\mathbb{E}\!\left[\|\nabla f(\bx_k)\|^2\right]
+ \frac{10T r^2L^2d^2}{\tilde{\beta}^2}
+ \frac{5d}{2\tilde{\beta}}\|\nabla f(\bx_1)\|^2
,
\end{align*}
After plugging in the bound~\eqref{eq:bound_pk}, we can show that
\begin{align*}
\left(1-\frac{4\eta^2 d^2G^2(1+\alpha^2)}{\tilde{\beta}^2 r^2(1-\alpha^2)^2}\right)\sum_{k=1}^T\mathbb{E}\!\left[\|\bg_k\|^2\right]
\leq\ &
\frac{5d}{\tilde{\beta}^2}\sum_{k=1}^{T}\mathbb{E}\!\left[\|\nabla f(\bx_k)\|^2\right]
+ \frac{10Tr^2L^2d^2}{\tilde{\beta}^2}
+\frac{5d}{2\tilde{\beta}}\|\nabla f(\bx_1)\|^2
,
\end{align*}
which gives~\eqref{eq:bound_gk}.
\end{proof}

\section{Proof of Theorem~\ref{theorem:main_convex}}\label{appendix:proof_convex}


Denote $\bw_k\coloneqq \bx_k - \alpha \bp_{k-1}/(1-\alpha)$. We have
$$
\bw_{k+1}
=\bx_{k+1} - \frac{\alpha}{1\!-\!\alpha}\bp_k =
\bx_k - \frac{1}{1\!-\!\alpha}\bp_k =
\bx_k - \frac{\alpha}{1\!-\!\alpha}\bp_{k-1} - \frac{\eta}{1\!-\!\alpha} \bg_k
=\bw_{k} - \frac{\eta}{1\!-\!\alpha} \bg_k,
$$
and thus
$$
\begin{aligned}
\left\|\bw_{k+1}-\bx^\ast\right\|^2
=\ &
\left\|\bw_k-\bx^\ast\right\|^2
-\frac{2\eta}{1-\alpha}\left\langle \bw_k-\bx^\ast, \bg_k\right\rangle
+\frac{\eta^2}{(1-\alpha)^2}\|\bg_k\|^2.
\end{aligned}
$$
By taking the expectation conditioned on $\mathcal{F}_k$, we get
\begin{align*}
\mathbb{E}\!\left[\!\left.\|\bw_{k+1}-\bx^\ast\|^2\,\right|\mathcal{F}_k\right]
=\ & \|\bw_k - \bx^\ast\|^2
-\frac{2\eta}{1\!-\!\alpha}\langle \bw_k-\bx^\ast, \nabla f_r(\bx_k)\rangle +
\frac{\eta^2}{(1\!-\!\alpha)^2}\mathbb{E}\!\left[\!
\left.\|\bg_k\|^2\,\right|\mathcal{F}_k
\right] \\
=\ & 
\|\bw_k - \bx^\ast\|^2
- \frac{2\eta}{1\!-\!\alpha}\langle \bx_k-\bx^\ast, \nabla f_r(\bx_k)\rangle \\
&
- \frac{2\eta\alpha}{(1\!-\!\alpha)^2}\langle \bx_k-\bx_{k-1},\nabla f_r(\bx_k)\rangle
+ \frac{\eta^2}{(1\!-\!\alpha)^2}\mathbb{E}\!\left[\!\left.\|\bg_k\|^2\,\right|\mathcal{F}_k
\right].
\end{align*}
Noticing that the convexity of $f$ implies
$$
-\langle \bx_k - \bx_{k-1},\nabla f_r(\bx_k)\rangle
\leq -(f_r(\bx_k)-f_r(\bx_{k-1})),
$$
$$
-\langle \bx_k-\bx^\ast,\nabla f_r(\bx_k)\rangle
\leq -(f_r(\bx_k)-f_r(\bx^\ast)),
$$
we get
\begin{align*}
\mathbb{E}\!\left[\|\bw_{k+1}-\bx^\ast\|^2\right]
\leq\ &
\mathbb{E}\!\left[\|\bw_{k}-\bx^\ast\|^2\right]
-\frac{2\eta}{1\!-\!\alpha}
\mathbb{E}\!\left[f_r(\bx_k) \!-\! f_r(\bx^\ast)\right] \\
&
-\frac{2\eta\alpha}{(1-\alpha)^2}
\mathbb{E}\!\left[f_r(\bx_k) - f_r(\bx_{k-1})\right]
+\frac{\eta^2}{(1\!-\!\alpha)^2}\mathbb{E}\!\left[\|\bg_k\|^2\right].
\end{align*}
Now let us assume $\theta \coloneqq 1-4\eta^2 d^2 G^2(1+\alpha^2)/[\tilde{\beta}^2 r^2(1-\alpha^2)^2]>0$. By taking the telescoping sum over $k=1,\ldots,T$ and plugging in the bound on $\sum_{k=1}^T \mathbb{E}\!\left[\|\bg_k\|^2\right]$ in Lemma~\ref{lemma:bound_gk_pk}, we get
\begin{align*}
\mathbb{E}\!\left[\|\bw_{T+1} \!-\! \bx^\ast\|^2\right]
\leq\ & \|\bw_{1} \!-\! \bx^\ast\|^2
-\frac{2\eta}{1 \!-\! \alpha}\sum_{k=1}^T
\mathbb{E}\!\left[f_r(\bx_k) \!-\! f_r(\bx^\ast)\right]
-\frac{2\eta\alpha}{(1 \!-\! \alpha)^2}\mathbb{E}[f_r(\bx_T) \!-\! f_r(\bx_1)] \\
& +\frac{10\eta^2Ld}{(1 \!-\! \alpha)^2\theta\tilde{\beta}^2}
\sum_{k=1}^T\mathbb{E}[f(\bx_k)-f(\bx^\ast)]
+ \frac{10T\eta^2 r^2L^2d^2}{(1 \!-\! \alpha)^2\theta\tilde{\beta}^2}
+\frac{5\eta^2Ld}{(1\!-\!\alpha)^2\theta\tilde{\beta}}(f(\bx_1)-f(\bx^\ast))
,
\end{align*}
where we also used the inequality $\|\nabla f(\bx_k)\|^2\leq 2L(f(\bx_k)-f(\bx^\ast))$ (see equation (2.1.7) in \cite{nesterov2004introductory}). By Lemma~\ref{lemma:function_fr}, we get
\begin{align*}
\mathbb{E}\!\left[\|\bw_{T+1} \!-\! \bx^\ast\|^2\right]
\leq\ & \|\bw_{1} \!-\! \bx^\ast\|^2
-\frac{2\eta}{1 \!-\! \alpha}\sum_{k=1}^T
\mathbb{E}\!\left[f(\bx_k) \!-\! f(\bx^\ast)\right]
-\frac{2\eta\alpha}{(1 \!-\! \alpha)^2}\mathbb{E}[f(\bx_T) \!-\! f(\bx_1)]
+ \frac{2\eta T r^2 L}{1 \!-\! \alpha}
\!+\! \frac{2\eta\alpha r^2 L}{(1 \!-\! \alpha)^2}
\\
&\!\!\!\!
+\frac{10\eta^2Ld}{(1 \!-\! \alpha)^2\theta\tilde{\beta}^2}\!
\sum_{k=1}^T\mathbb{E}[f(\bx_k) \!-\! f(\bx^\ast)]
+ \frac{10T\eta^2 r^2L^2d^2}{(1 \!-\! \alpha)^2\theta\tilde{\beta}^2}
+\frac{5\eta^2 Ld}{(1\!-\!\alpha)^2\theta\tilde{\beta}}
(f(\bx_1)-f(\bx^\ast))
,
\end{align*}
which further leads to
\begin{align*}
& \left(
1-\frac{5\eta Ld}{(1\!-\!\alpha)\theta\tilde{\beta}^2}
\right)\frac{1}{T}\sum_{k=1}^T
\mathbb{E}[f(\bx_k)-f(\bx^\ast)]\\
\leq \ & \left(
1-\frac{5\eta Ld}{(1\!-\!\alpha)\theta\tilde{\beta}^2}
\right)\frac{1}{T}\sum_{k=1}^T
\mathbb{E}[f(\bx_k)-f(\bx^\ast)]
+\frac{\alpha}{1\!-\!\alpha}\frac{\mathbb{E}[f(\bx_T)-f(\bx^\ast)]}{T} \\
\leq\ &
\frac{(1\!-\!\alpha)\|\bx_1-\bx^\ast\|^2}{2\eta T}
+\frac{5\eta r^2L^2d^2}{(1\!-\!\alpha)\theta\tilde{\beta}^2}
+ r^2L  + \frac{\alpha r^2 L}{(1\!-\!\alpha)T}
+\left(\frac{\alpha}{1\!-\!\alpha}
+\frac{5\eta Ld}{2(1\!-\!\alpha)\theta\tilde{\beta}}
\right)\frac{f(\bx_1)-f(\bx^\ast)}{T}.
\end{align*}
Now, by letting
$$
\eta \leq \frac{(1-\alpha)\tilde{\beta}^2}{20Ld T^{1/3}},\qquad
\frac{4\eta d G}{\tilde{\beta}(1\!-\!\alpha)}\leq r \leq \frac{G}{L T^{1/3}},
$$
we see that
$$
\theta \geq 1 - \frac{1}{4}
\cdot \frac{(1+\alpha^2)}{(1+\alpha)^2}\geq \frac{3}{4},
\qquad
1-\frac{5\eta Ld}{(1\!-\!\alpha)\theta\tilde{\beta}^2}
\geq 1-\frac{1}{3T^{1/3}}\geq \frac{2}{3},
$$
and then we obtain
\begin{align*}
& \frac{1}{T}\sum_{k=1}^T
\mathbb{E}[f(\bx_k)-f(\bx^\ast)] \\
\leq\ &
\frac{3(1\!-\!\alpha)\|\bx_1-\bx^\ast\|^2}{4\eta T} + \frac{G^2 d}{2 LT}
+\frac{3G^2}{2 L T^{2/3}} + \frac{3\alpha  G^2}{2(1-\alpha) L T^{5/3}}
+ \frac{3}{2}\left(\frac{\alpha}{1\!-\!\alpha}
+\frac{10\eta Ld}{3(1\!-\!\alpha)\tilde{\beta}}
\right)\frac{f(\bx_1)-f(\bx^\ast)}{T},
\\
\leq\ &
\frac{3(1\!-\!\alpha)\|\bx_1-\bx^\ast\|^2}{4\eta T}
+\frac{3G^2}{2 L T^{2/3}}
+ \sO(d/T),
\end{align*}
which implies the desired bound as $f$ is convex.

\section{Proof of Theorem~\ref{theorem:main_nonconvex}}\label{appendix:proof_nonconvex}


We  let $\bw_k\coloneqq \bx_k-\alpha \bp_{k-1}/(1-\alpha)$, and it can be checked that $\bw_{k+1}=\bw_k-\eta \bg_k/(1-\alpha)$. By the $L$-smoothness of $f$, we have
\begin{align*}
f(\bw_{k+1})\leq\ &
f(\bw_k) -
\frac{\eta}{1 \!-\! \alpha}\langle \bg_k,\nabla f(\bw_k)\rangle
+\frac{\eta^2 L}{2(1 \!-\! \alpha)^2}\|\bg_k\|^2 \\
=\ &
f(\bw_k) -
\frac{\eta}{1 \!-\! \alpha}\langle \bg_k,\nabla f(\bw_k)-\nabla f(\bx_k)\rangle
-\frac{\eta}{1 \!-\! \alpha}\langle \bg_k, \nabla f(\bx_k)\rangle
+\frac{\eta^2 L}{2(1 \!-\! \alpha)^2}\|\bg_k\|^2.
\end{align*}
By taking the expectation conditioned on $\mathcal{F}_k$, we get
\begin{align*}
\mathbb{E}[f(\bw_{k+1})|\mathcal{F}_k]
\leq\ &
f(\bw_k)
- \frac{\eta}{1 \!-\! \alpha}\langle \nabla f_r(\bx_k),\nabla f(\bw_k) \!-\! \nabla f(\bx_k)\rangle
-\frac{\eta}{1 \!-\! \alpha}\langle\nabla f_r(\bx_k),\nabla f(\bx_k)\rangle \\
& + \frac{\eta^2L}{2(1 \!-\! \alpha)^2}
\mathbb{E}[\|\bg_k\|^2|\mathcal{F}_k].
\end{align*}
Note that
\begin{align*}
& -\langle \nabla f_r(\bx_k), \nabla f(\bw_k)-\nabla f(\bx_k)\rangle \\
=\ & -\langle \nabla f(\bx_k), \nabla f(\bw_k)-\nabla f(\bx_k)\rangle
-\langle \nabla f_r(\bx_k)-\nabla f(\bx_k), \nabla f(\bw_k)-\nabla f(\bx_k)\rangle \\
\leq \ & \frac{1}{4}\|\nabla f(\bx_k)\|^2 + \|\nabla f(\bw_k) - \nabla f(\bx_k)\|^2 + \|\nabla f_r(\bx_k) - \nabla f(\bx_k)\|^2 + \frac{1}{4}\|\nabla f(\bw_k) - \nabla f(\bx_k)\|^2\\
=\ &
\frac{1}{4}\|\nabla f(\bx_k)\|^2
+ \|\nabla f_r(\bx_k) - \nabla f(\bx_k)\|^2
+ \frac{5}{4}\|\nabla f(\bw_k) - \nabla f(\bx_k)\|^2 \\
\leq\ &
\frac{1}{4} \|\nabla f(\bx_k)\|^2
+ L^2r^2
+ \frac{5}{4}L^2 \|\bw_k-\bx_k\|^2   \tag{by Lemma \ref{lemma:function_fr} and $L$-Smoothness of $f$} \\
=\ &
\frac{1}{4} \|\nabla f(\bx_k)\|^2
+ L^2r^2
+ \frac{5L^2\alpha^2}{4(1\!-\!\alpha)^2}\|\bp_{k-1}\|^2,
\end{align*}
and
\begin{align*}
-\langle\nabla f_r(\bx_k),\nabla f(\bx_k)\rangle
=\ &
-\|\nabla f(\bx_k)\|^2
-\langle\nabla f_r(\bx_k) - \nabla f(\bx_k),\nabla f(\bx_k)\rangle \\
\leq \ &   - \|\nabla f(\bx_k)\|^2 + \|\nabla f_r(\bx_k)-\nabla f(\bx_k)\|^2 +  \frac{1}{4}\|\nabla f(\bx_k)\|^2 \\
=\ &
-\frac{3}{4}\|\nabla f(\bx_k)\|^2
+\|\nabla f_r(\bx_k)-\nabla f(\bx_k)\|^2 \\
\leq\ &
-\frac{3}{4}\|\nabla f(\bx_k)\|^2
+ L^2r^2.
\end{align*}
Therefore
\begin{align*}
\mathbb{E}\!\left[f(\bw_{k+1})|\mathcal{F}_k\right]
\leq\ &
f(\bw_{k})
-\frac{\eta}{2(1\!-\!\alpha)}\|\nabla f(\bx_k)\|^2
+\frac{2\eta r^2 L^2}{1\!-\!\alpha} \\
& +\frac{5\eta L^2\alpha^2}{4(1\!-\!\alpha)^3}\|\bp_{k-1}\|^2
+\frac{\eta^2 L}{2(1\!-\!\alpha)^2}
\mathbb{E}[\|\bg_k\|^2|\mathcal{F}_k].
\end{align*}
Now, by taking the total expectation and summing over $k=1,\ldots,T$, we get
\begin{align*}
\mathbb{E}[f(\bw_{T+1})]
\leq\ &
f(\bw_1) - \frac{\eta}{2(1\!-\!\alpha)}\sum_{k=1}^T \mathbb{E}[\|\nabla f(\bx_k)\|^2]
+\frac{2\eta r^2 L^2 T}{1\!-\!\alpha} \\
&
+\frac{5\eta L^2\alpha^2}{4(1\!-\!\alpha)^3}
\sum_{k=1}^{T-1}\mathbb{E}[\|\bp_{k}\|^2]
+\frac{\eta^2 L}{2(1\!-\!\alpha)^2}
\sum_{k=1}^T\mathbb{E}[\|\bg_k\|^2].
\end{align*}
Assuming $\theta\coloneqq 1-4\eta^2 d^2G^2(1+\alpha^2)/[\tilde{\beta}^2 r^2(1-\alpha^2)^2] >0$ and plugging in the bounds in Lemma~\ref{lemma:bound_gk_pk}, we get
\begin{align*}
\mathbb{E}[f(\bw_{T+1})]
\leq\ &
f(\bw_1) - \frac{\eta}{2(1\!-\!\alpha)}\sum_{k=1}^T \mathbb{E}[\|\nabla f(\bx_k)\|^2]
+\frac{2\eta r^2 L^2 T}{1\!-\!\alpha} \\
&
+\frac{\eta^2 L}{2(1\!-\!\alpha)^2}
\!\left(\! 1 \!+ \!
\frac{5\eta L\alpha^2(1\!+\!\alpha^2)}{(1\!-\!\alpha)^3(1\!+\!\alpha)^2}
\right)
\!\!\left(\frac{5d}{\theta\tilde{\beta}^2}\sum_{k=1}^T \mathbb{E}\!\left[\|\nabla f(\bx_k)\|^2\right]
+ \frac{10 T r^2 L^2 d^2}{\theta\tilde{\beta}^2}
+\frac{5d}{2\theta\tilde{\beta}}\|\nabla f(\bx_1)\|^2 
\right).
\end{align*}
Now let $T$ be sufficiently large such that
$$
T^{1/3}\geq \frac{\alpha^2 \tilde{\beta}^2}{2d(1-\alpha)^2},
$$
and let
$$
\eta \leq \frac{(1\!-\!\alpha)\tilde{\beta}^2}{20Ld T^{1/3}},\qquad
\frac{4\eta d G}{\tilde{\beta}(1\!-\!\alpha)}\leq r \leq \frac{G}{L T^{1/3}},
$$
We then have $\eta \leq (1-\alpha)^3/(10\alpha^2 L)$, and thus
$$
\theta \geq 1-\frac{1}{4}\cdot\frac{(1\!+\!\alpha^2)}{(1\!+\!\alpha)^2}\geq\frac{3}{4},
\qquad
\frac{5\eta L\alpha^2(1\!+\!\alpha^2)}{(1\!-\!\alpha)^3(1\!+\!\alpha)^2}
\leq\frac{1\!+\!\alpha^2}{2(1\!+\!\alpha)^2}
\leq \frac{1}{2}.
$$
Therefore
\begin{align*}
\mathbb{E}[f(\bw_{T+1})]
\leq\ &
f(\bw_1) - \frac{\eta}{2(1\!-\!\alpha)}
\left(1-\frac{10\eta Ld}{(1\!-\!\alpha)\tilde{\beta}^2}\right)\sum_{k=1}^T \mathbb{E}[\|\nabla f(\bx_k)\|^2] \\
&
+\frac{2\eta r^2 L^2 T}{1\!-\!\alpha}
+
\frac{\eta^2 L}{(1\!-\!\alpha)^2\tilde{\beta}^2}\cdot 10 T r^2 L^2 d^2
+\frac{5\eta^2 Ld}{2(1\!-\!\alpha)^2\tilde{\beta}}\|\nabla f(\bx_1)\|^2
\\
\leq\ &
f(\bx_1) - \frac{\eta}{4(1\!-\!\alpha)}
\sum_{k=1}^T \mathbb{E}[\|\nabla f(\bx_k)\|^2] 
+\frac{\eta}{1\!-\!\alpha}\cdot 2 G^2 T^{1/3}
+
\frac{\eta}{1\!-\!\alpha}
\cdot \frac{G^2d}{2}
+\frac{\eta}{1\!-\!\alpha}\frac{\tilde{\beta}}{8T^{1/3}}\|\nabla f(\bx_1)\|^2.
\end{align*}
Since $f^\ast\leq \mathbb{E}[f(\bw_{T+1})]$, 
it leads to
\begin{align*}
\frac{1}{T}\sum_{k=1}^T
\mathbb{E}\!\left[
\|\nabla f(\bx_k)\|^2
\right]
\leq\ &
\frac{4(1\!-\!\alpha)(f(\bx_1)-f^\ast)}{\eta T}
+\frac{8 G^2}{T^{2/3}}
+\frac{2G^2 d}{ T}
+\frac{\tilde{\beta}\|\nabla f(\bx_1)\|^2}{2T^{4/3}}
.
\end{align*}

\end{document}